\documentclass{l4dc2024}
\usepackage{natbib}
\usepackage{multirow}
\usepackage{mathbbol}
\usepackage{array}
\usepackage{amssymb,arydshln}
\usepackage[nodayofweek]{datetime}
\usepackage{float}
\usepackage{bm}
\usepackage{threeparttable}
\usepackage{booktabs}
\usepackage{balance}
\usepackage{tcolorbox}
\usepackage{xcolor}
\usepackage[font=small,labelfont=bf]{caption} 
\usepackage{graphicx,color}
\usepackage{subcaption,setspace}
\usepackage{multirow}
\usepackage{wrapfig}
\usepackage{optidef}

\usepackage{algorithm}
\usepackage{algorithmicx}
\usepackage{algpseudocode}
\usepackage{mathrsfs}

\newcommand{\ini}{\textnormal{ini}}
\newcommand{\col}{\textnormal{col}}
\newcommand{\f}{\textnormal{F}}
\newcommand{\p}{\textnormal{P}}
\newcommand{\D}{\textnormal{d}}
\newcommand{\R}{\textnormal{r}}
\newcommand{\h}{\textnormal{H}}
\newcommand{\tr}{{{\mathsf T}}}
\newcommand{\A}{\textnormal{a}}
\newcommand{\B}{\textnormal{b}}
\newcommand{\method}[1]{\texttt{#1}}

\newcommand{\DeePC}{\texttt{DeePC}}
\newcommand{\DeePCsvd}{\texttt{DeePC-SVD}}

\newtheorem{thm}{Theorem}
\newtheorem{prop}{Proposition}
\newtheorem{lem}{Lemma}

\newtheorem{fact}{Fact}

\usepackage[nameinlink]{cleveref}
\crefname{equation}{}{}
\crefname{theorem}{Theorem}{Theorems}
\crefname{corollary}{Corollary}{Corollaries}
\crefname{example}{Example}{Examples}
\crefname{assumption}{Assumption}{Assumptions}
\crefname{lemma}{Lemma}{Lemmas}
\crefname{proposition}{Proposition}{Propositions}
\crefname{figure}{Fig.}{Figures}
\crefname{table}{Table}{Tables}
\crefname{section}{Section}{Sections}
\crefname{appendix}{Appendix}{Appendices}
\crefname{definition}{Definition}{Definitions}

\crefname{thm}{Theorem}{Theorems}
\crefname{prop}{Proposition}{Propositions}
\crefname{lem}{Lemma}{Lemmas}
\crefname{deff}{Definition}{Definitions}
\crefname{rem}{Remark}{Remarks}
\crefname{col}{Corollary}{Corollaryies}

\Crefname{equation}{}{}
\Crefname{theorem}{Theorem}{Theorems}
\Crefname{corollary}{Corollary}{Corollaries}
\Crefname{example}{Example}{Examples}
\Crefname{lemma}{Lemma}{Lemmas}
\Crefname{proposition}{Proposition}{Propositions}
\Crefname{figure}{Figure}{Figures}
\Crefname{fact}{Fact}{Facts}
\Crefname{table}{Table}{Tables}
\Crefname{section}{Section}{Sections}
\Crefname{appendix}{Appendix}{Appendices}
\crefname{prop}{Proposition}{Propositions}

\title[Convex Approximations for DeePC]{Convex Approximations for a Bi-level Formulation of Data-Enabled Predictive Control}
\usepackage{times}

\author{%
 \Name{Xu Shang} \Email{x3shang@ucsd.edu}\\
 \Name{Yang Zheng} \Email{zhengy@ucsd.edu}\\
 \addr Department of Electrical and Computer Engineering, University of California San Diego%
}

\begin{document}

\maketitle

\vspace{-8mm}
\begin{abstract}%

The Willems' fundamental lemma, which characterizes linear time-invariant (LTI) systems using input and output trajectories, has found many successful applications. Combining this with receding horizon control leads to a popular Data-EnablEd Predictive Control (\DeePC{}) scheme. \DeePC{} is first established for LTI systems and has been extended and applied for practical systems beyond LTI settings. However, the relationship between different \DeePC{} variants, involving regularization and dimension reduction, remains unclear. In this paper, we first introduce a new bi-level optimization formulation that combines a data pre-processing step as an inner problem (system identification) and predictive control as an outer problem (online control). We next discuss a series of convex approximations by relaxing some hard constraints in the bi-level optimization as suitable regularization terms, accounting for an implicit identification. These include some existing \DeePC{} variants as well as two new variants, for which we establish their equivalence under appropriate settings. 
Notably, our analysis reveals a novel variant, called \method{DeePC-SVD-Iter}, which has remarkable empirical performance of direct methods on systems beyond deterministic LTI settings.     
\end{abstract}


\section{Introduction}
There has been a surging interest in utilizing data-driven techniques to control systems with unknown dynamics \citep{pillonetto2014kernel, markovsky2021behavioral}. Existing methods can be generally categorized into indirect and direct data-driven control techniques: indirect data-driven control approaches typically include sequential system identification (system ID) \citep{ljung1998system, chiuso2019system} and model-based control \citep{kouvaritakis2016model}, while direct data-driven control methods bypass system ID and directly design control strategies from input and output measured data \citep{markovsky2021behavioral}. 

In particular, Data-EnablEd Predictive Control (\DeePC{}{}) \citep{coulson2019data, markovsky2021behavioral} that combines behavioral theory with receding horizon control has received increasing attention. It utilizes Willem's fundamental lemma \citep{willems2007behavior} to construct a data-driven representation of a dynamic system and incorporates it with receding horizon control. \DeePC{}{} is first established for deterministic linear time-invariant (LTI) systems, and its equivalence with subspace predictive control (SPC) has been discussed in \citep{fiedler2021relationship}. \cite{berberich2020data} further investigate conditions for its closed-loop stability. The \DeePC{} approach has shown promising results for the control of practical systems beyond LTI settings \citep{wang2023deep, elokda2021data, shang2023smoothing,lian2023adaptive}. For non-deterministic or nonlinear systems, suitable regularizations are necessary for \DeePC{}; see \cite{breschi2023data, dorfler2022bridging}.

\begin{figure}[t]
    \centering
    \includegraphics[width = \textwidth]{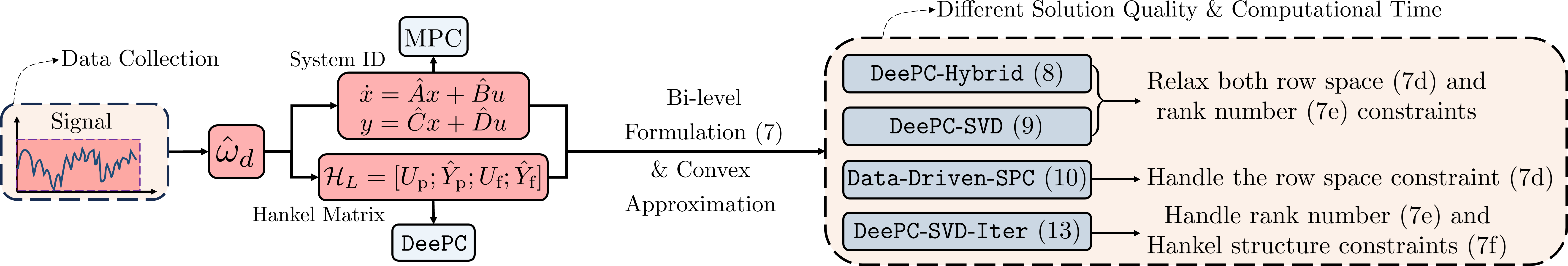}
    \caption{Schematic of data-driven control, which starts by collecting data (usually noisy) from the real system. Indirect methods identify a parametric model, while \DeePC{} forms a Hankel matrix as the trajectory library for predictive control. Our bi-level formulation \cref{eqn:bi-level} integrates system ID techniques to \DeePC{}. We introduce a series of convex approximations \cref{eqn:DeePC-Hybrid}, \cref{eqn:DeePC-SVD}, \cref{eqn:data-driven-spc}, and \cref{eqn:DeePC-SVD-Iter} that relax the bi-level formulation.}
    \label{fig:WorkFlow}
    \vspace{-6mm}
\end{figure}

There are different regularization strategies for \DeePC{}, ranging from some heuristics in \cite{coulson2019data} to principled analysis via bi-level formulations in \cite{dorfler2022bridging}. Notably, indirect data-driven control is first formulated as a bi-level optimization problem involving both control and identification in \cite{dorfler2022bridging}. Many regularized versions (such as $l_1$ or $l_2$~norms) of \DeePC{} can be considered as convex relaxations of the bi-level optimization. Beyond regularization, some recent approaches aim to decrease the optimization dimensions in \DeePC{} and improve computational efficiency \citep{zhang2023dimension, alsalti2023data}. One simple strategy is to use a singular value decomposition (SVD) to pre-process the data-driven representation, which has shown promising performance \citep{zhang2023dimension}. However, the relationship between the recent variants of \DeePC{} for non-deterministic and nonlinear systems, involving regularization and dimension reduction, remains unclear, and there is no analysis and comparison for their solution qualities.

In this paper, we introduce a new bi-level formulation incorporating both system ID techniques and predictive control, and discuss how existing and new variants of \DeePC{} can be considered as convex approximations of this bi-level formulation; \Cref{fig:WorkFlow} illustrates the overall process. Specifically, in our bi-level formulation for \DeePC{}, the data pre-processing step is viewed as an inner optimization problem (identification), and the predictive control is viewed as an outer optimization problem (online control). Constraints for the inner optimization problem are derived from system ID methods (e.g., SPC \cite{favoreel1999spc}, and low-rank approximation \cite{markovsky2016missing}). We further discuss a series of convex approximations by relaxing some hard constraints as suitable regularization terms. In this process, we derive two new variants of \DeePC{} by adapting existing methods: 1) \method{Data-Driven-SPC} is derived from classical SPC with the same structure as \DeePC{}, and 2) \method{DeePC-SVD-Iter} refines the data-driven representation in \method{DeePC-SVD} from \cite{zhang2023dimension} and provides superior performance. We also investigate the equivalence of \method{DeePC-Hybrid} \citep{dorfler2022bridging}, \DeePCsvd{} \citep{zhang2023dimension}, and \method{Data-Driven-SPC}. Our analysis is more general than \cite{zhang2023dimension,fiedler2021relationship,breschi2023data}. Numerical experiments confirm our analysis and show the superior performance of \method{DeePC-SVD-Iter}. 

The rest of this paper is structured as follows. \Cref{sec:Prob-Meth} introduces preliminaries and the problem statement of our bi-level formulation. A series of convex approximations are introduced in \Cref{sec:methods}, their relationship is established in \Cref{sec:relationships}. \Cref{sec:results} compares their control performance via numerical simulations. Finally, we conclude the paper in \Cref{sec:conclusions}. Some notations and background are listed in \Cref{appendix:review}. We use $\col(A_1, A_2, \ldots, A_m) = [
A_1^\tr, A_2^\tr, \ldots, A_m^\tr]^\tr
$.

\section{Preliminaries and Problem Statement}
\label{sec:Prob-Meth}
\subsection{Preliminaries}
We consider a linear time-invariant (LTI) system in the discrete-time domain:
\begin{equation}
\label{eqn:LTI-sys}
    \left\{\begin{aligned}
    x(k+1) &= A x(k) + B u(k), \\
    y(k) &= C x(k) + D u(k),
    \end{aligned} \right.
\end{equation}
where the state, input, output at time $k$ are $x(k) \in \mathbb{R}^n$, $u(k) \in \mathbb{R}^m$, and $y(k) \in \mathbb{R}^p$, respectively. 
Given a desired reference trajectory $y_\textrm{r} \in \mathbb{R}^{pN}$ with horizon $N >0$, input constraint set $\mathcal{U} \subseteq \mathbb{R}^m$, output constraint set $
\mathcal{Y} \subseteq \mathbb{R}^p$, we aim to design control inputs such that the system output tracks the reference trajectory. In particular, we consider the well-known receding horizon predictive control
\begin{subequations}
\label{eqn:predictive-control}
{
\begin{align}
    \min_{x, u ,y } \quad & \sum_{k=t}^{t+N-1}\big(\|y(k) - y_\textrm{r}(k)\|_Q^2 + \|u(k)\|_R^2\big) \nonumber \\ 
    \textrm{subject~to} \quad & x(k+1) = A x(k)+ B u(k), \quad k \in [t, t+N-1] \label{eqn:predictive-control-a} \\
    & y(k) = C x(k) + D u(k), \quad k \in [t, t+N-1] \label{eqn:predictive-control-b} \\
    & x(t) = x_\ini, \label{eqn:predictive-control-c}\\
    & u(k) \in \mathcal{U}, \; y(k) \in \mathcal{Y}, \quad k \in [t, t+N-1],  \label{eqn:predictive-control-d}
\end{align}}
\end{subequations}

\noindent where $x_\ini \! \in \! \mathbb{R}^n$ is the initial state of system \cref{eqn:LTI-sys} and $\|u(k)\|_R^2$ denotes the quadratic norm $u(k)^\tr \! R u(k)$ (similarly for $\|\cdot\|_Q$) with $R\in \mathbb{S}^m_+$ and $Q \in \mathbb{S}^p_+$. 
We assume $\mathcal{U}$ and $\mathcal{Y}$ are convex sets. Without loss of generality, we consider a regulation problem (i.e., $y_\textrm{r} = \mathbb{0}_{pN}$) from the rest of the discussions. 

It is clear that \cref{eqn:predictive-control} is a convex optimization problem (it is indeed a quadratic program for simple $\mathcal{U}$ and $\mathcal{Y}$), which admits an efficient solution when the model for the system \cref{eqn:LTI-sys} is known, i.e.,
matrices $A$, $B$, $C$ and $D$ are known. In this work, we focus on the case when the system model and the initial condition $x_\ini$ are unknown. Instead, we have access to 1) offline data, i.e., a length-$T$ pre-collected input/output trajectory of \cref{eqn:LTI-sys}, and 2) online data, i.e., the most recent past input/output sequence of length-$T_\ini$. Then, \cref{eqn:predictive-control} can be solved  by either indirect system ID and model-based control \citep{aastrom1971system} 
or the recent emerging direct data-driven control, such as \DeePC{} and its related approaches \citep{dorfler2022bridging,markovsky2021behavioral}. As discussed in \cite{dorfler2022bridging}, the indirect system ID approach is superior in the case of ``variance'' noise, while \DeePC{} with suitable regularization terms has better performance in the case of ``bias'' errors. 

\subsection{Data-Enabled Predictive Control}
\label{subsec:DeePC}
We here review the basic setup of \DeePC{}. First, let us introduce a notion of persistent excitation. 
\begin{definition}[Persistently Exciting]
\label{def:Hankel-struct}
    A sequence of signal $\omega  =  \textnormal{col}(\omega(1),\omega(2), \ldots, \omega(T))$ of the length $T$ ($T \in \mathbb{N}$) is persistently exciting of order $L$ ($L < T$) if its associated Hankel matrix with depth $L$,  defined below, has full row rank, 
    {\small \[\mathcal{H}_L(\omega) = \begin{bmatrix}
        \omega(1) & \omega(2) & \cdots & \omega(T-L+1) \\
        \omega(2) & \omega(3) & \cdots & \omega(T-L+2) \\
        \vdots    & \vdots    & \ddots & \vdots \\
        \omega(L) & \omega(L+1) & \cdots & \omega(T)
    \end{bmatrix}. \]} \vspace{-5mm}
\end{definition}

\begin{lem}[Fundamental Lemma; \cite{willems2005note}] \label{lemma:Fundamental}
    Suppose that system \cref{eqn:LTI-sys} is controllable. Given a length $T$ input/output trajectory: $u_\textnormal{d} = \textnormal{col}(u_\textnormal{d}(1),\ldots, u_\textnormal{d}(T)) \in \mathbb{R}^{mT}$, $
    y_\textnormal{d} = \textnormal{col}(y_\textnormal{d}(1)$, $ \ldots, y_\textnormal{d}(T)) \in \mathbb{R}^{pT}$ where $u_\textnormal{d}$ is persistently exciting of order $L+n$, then a length $L$ input/output sequence $(u_\textnormal{s}, y_\textnormal{s})$ is a valid trajectory of \cref{eqn:LTI-sys} if and only if there exists a $g \in \mathbb{R}^{T-L+1}$ such that 
    {\small  \begin{equation} 
    \label{eqn:FundaLemma}
    \begin{bmatrix}
    \mathcal{H}_L(u_\textnormal{d})\\
    \mathcal{H}_L(y_\textnormal{d})
    \end{bmatrix} g
    =
    \begin{bmatrix}
    u_\textnormal{s}\\
    y_\textnormal{s}
    \end{bmatrix}. \vspace{-1mm}
    \end{equation}}
    If length $L$ is not smaller than the lag of the system, matrix $\textnormal{col}(\mathcal{H}_L(u_\textnormal{d}), \mathcal{H}_L(y_\textnormal{d}))$ has rank $mL+n$.
\end{lem}

The \DeePC{} approach in \cite{coulson2019data} employs \cref{eqn:FundaLemma} to build a predictor based on the pre-collected data. 
In particular, the Hankel matrix formed by the offline data is partitioned as
 \begin{equation}
\label{eqn:Hankel_partition}
  {\small      \begin{bmatrix}
    U_{\textrm{P}} \\
    U_{\textrm{F}} 
\end{bmatrix} := \mathcal{H}_L(u_\textrm{d}), \quad 
\begin{bmatrix}
    Y_{\textrm{P}} \\
    Y_{\textrm{F}} 
\end{bmatrix} := \mathcal{H}_L(y_\textrm{d}),}
\end{equation}
where $U_\p$ and $U_\f$ consist the first $T_\ini$ rows and the last $N$ rows of $\mathcal{H}_L(u_\D)$, respectively (similarly for $Y_\p$ and $Y_\f$; so $L = T_\ini + N$). We denote the most recent past input trajectory of length $T_\ini$ and the future input trajectory of length $N$ as $u_\ini = \col(u(t-T_\textrm{ini}),u(t-T_\textrm{ini}+1),\ldots, u(t-1))$ and $u = \col(u(t), u(t+1),\ldots, u(t+N-1))$, respectively (similarly for $y_\textrm{ini}, y$).

\begin{wrapfigure}[4]{r}{.3\textwidth}
\vspace{-4.5mm}
{\small \begin{equation}
\label{eqn:predictor}
\begin{bmatrix}
    U_\textrm{P} \\
    Y_\textrm{P} \\
    U_\textrm{F} \\
    Y_\textrm{F}
\end{bmatrix} g
=
\begin{bmatrix}
    u_\textrm{ini} \\
    y_\textrm{ini} \\
    u \\
    y
\end{bmatrix}.
\end{equation}}
\end{wrapfigure}
Then, \Cref{lemma:Fundamental} ensures that the sequence $\col(u_\textrm{ini}, y_\textrm{ini}, u, y)$ is a valid trajectory of \cref{eqn:LTI-sys} if and only if there exists $g \in \mathbb{R}^{T-T_\textrm{ini} -N+1}$ such that \cref{eqn:predictor} holds. 
For notational simplicity, we further denote the matrix $\col(U_\p, Y_\p, U_\f, Y_\f)$ associated with pre-collected data as $H$. Note that $H$ can be considered as a trajectory library since each of its columns is a valid trajectory of system \cref{eqn:LTI-sys}. If $T_{\mathrm{ini}}$ is larger or equal to the lag of system \cref{eqn:LTI-sys}, $y$ is unique given any $u_\textrm{ini}$, $y_\textrm{ini}$ and $u$ in \cref{eqn:predictor}. The most basic version of  \DeePC{} \citep{coulson2019data} utilizes the predictor \cref{eqn:predictor} as the data-driven representation of \cref{eqn:predictive-control-a,eqn:predictive-control-b,eqn:predictive-control-c} and reformulate the problem \cref{eqn:predictive-control} as
\begin{equation}
\label{eqn:DeePC}
\begin{aligned}
\min_{g,u,y} \quad & \sum_{k=t}^{t+N-1}\big(\|y(k) - y_r(k)\|_Q^2 + \|u(k)\|_R^2\big) \\
\textrm{subject to} \quad & \cref{eqn:predictor}, u \in \mathcal{U} ,y \in \mathcal{Y}
\end{aligned}
\end{equation} 
where we slightly abuse the notation and use $u \in \mathcal{U} ,y \in \mathcal{Y}$ to denote input/output constraints \cref{eqn:predictive-control-d}.

\subsection{A bi-level formulation beyond deterministic LTI systems}
It is not difficult to show that for LTI systems with noise-free data, problems \cref{eqn:predictive-control} and \cref{eqn:DeePC} are fully equivalent (cf. \citealp[Theorem 5.1]{coulson2019data}). However, for the case beyond deterministic~LTI systems, there exist different regularization terms or data pre-processing techniques that extend the basic \DeePC{} \cref{eqn:DeePC}. Indeed, an extensive discussion on bridging indirect and direct data-driven control was presented in \cite{dorfler2022bridging}, where two different bi-level formulations were discussed. 

Motivated by \cite{dorfler2022bridging}, we propose a new bi-level formulation that incorporates data pre-processing techniques from system ID. In practice, the data predictor $H$ in \cref{eqn:predictor} may be corrupted by  ``variance'' noises and/or ``bias'' errors. The key idea of our bi-level formulation is to pre-process the raw data $H$ and construct a new trajectory library $\Tilde{H}$ satisfying specific structures in system ID: 
\begin{mini!}[2]<b> 
        {g,\sigma_y,u \in \mathcal{U} ,y \in \mathcal{Y}}{\sum_{k=t}^{t+N-1}\big(\|y(k) - y_r(k)\|_Q^2 + \|u(k)\|_R^2\big) + \lambda_y \|\sigma_y\|_2^2}
        {\label{eqn:bi-level}}{\label{eqn:bi-level-a}}
        \addConstraint{\tilde{H}^* g}{=\col(u_\ini, y_\ini+\sigma_y, u, y ) \label{eqn:bi-level-b}}
        \addConstraint{\textrm{where}}{\;\; \tilde{H}^* \in \arg \min_{\tilde{H}} \quad J(\tilde{H},H) \label{eqn:bi-level-c}}
        \addConstraint{\phantom{\textrm{where}\;\; \tilde{H}^*} \mathrm{subject~to}}{\quad \tilde{Y}_\f = Y_\f /\col(\tilde{U}_\p, \tilde{Y}_\p, \tilde{U}_\f) \label{eqn:bi-level-d}}{ \ \text{(Row Space)}}
        \addConstraint{\phantom{\textrm{where}\;\; \tilde{H}^*  \mathrm{subject~to}}}{\quad \textrm{rank}(\tilde{H}) = mL+n \label{eqn:bi-level-e}}{\ \text{(Rank Number)}}
        \addConstraint{\phantom{\textrm{where}\;\; \tilde{H}^*  \mathrm{subject~to}}}{\quad \tilde{H} \in \mathcal{H} \label{eqn:bi-level-f}}{\ \text{(Hankel Structure).}} 
    \end{mini!}
This bi-level problem structure in \cref{eqn:bi-level}, which is consistent with those in \cite{dorfler2022bridging}, reflects the sequential
ID and control tasks, where we first fit a model $\tilde{H}$ from the raw data $H$ \cref{eqn:Hankel_partition} in the inner system ID before using the model for \DeePC{} in the outer problem.  

In the outer problem \cref{eqn:bi-level-a}-\cref{eqn:bi-level-b}, we have introduced a slack variable $\sigma_y$ and its regularization term to handle the model mismatch and ensure feasibility, as discussed in \cite{markovsky2021behavioral}. In the inner optimization problem \cref{eqn:bi-level-c}-\cref{eqn:bi-level-f}, $J(\tilde{H}, H)$ in \cref{eqn:bi-level-c} denotes system identification loss function with $H$ being the raw Hankel matrix \cref{eqn:Hankel_partition}. We have also partitioned the variable $\tilde{H}$ as $\col(\tilde{U}_\p, \tilde{Y}_\p, \tilde{U}_\f, \tilde{Y}_\f)$. In \cref{eqn:bi-level-d}, $Y_\f /\col(\tilde{U}_\p, \tilde{Y}_\p, \tilde{U}_\f)$ denotes the orthogonal projection of $Y_\f$ onto the row space of $\col(\tilde{U}_\p, \tilde{Y}_\p, \tilde{U}_\f)$. This row space constraint is derived from SPC \citep{favoreel1999spc} which will be discussed in detail in \Cref{subsec:data-driven-spc}. The rank constraint \cref{eqn:bi-level-e} and Hankel structure \cref{eqn:bi-level-f} (where $\mathcal{H}$ is the set of all matrices with Hankel structure; cf. Definition \ref{def:Hankel-struct}) come from low-rank approximation in \cite{markovsky2016missing}. We refer the interested readers to \cite{fiedler2021relationship} and \cite{willems2005note} for further details on row space and rank number respectively.   

We will derive a series of convex approximations for this bi-level formulation \cref{eqn:bi-level} in \Cref{sec:methods}  and discuss their equivalence (if possible) and relationship in \Cref{sec:relationships}. 

\section{Convex Approximations}
\label{sec:methods}
While the bi-level formulation \cref{eqn:bi-level} is not solvable immediately, it provides useful guidance~to~derive new formulations/variants of \DeePC{}. In this section, we present four convex approximations by adapting existing methods; see \Cref{fig:WorkFlow} for an overview. These strategies relax the inner constraints \cref{eqn:bi-level-d,eqn:bi-level-e,eqn:bi-level-f} using suitable regularizers to the outer problem.  
  
\subsection{\DeePC{} with regularization and dimension reduction}
\label{subsec:DeePC-reg}
We first discuss two existing convex approximations of \cref{eqn:bi-level}: \method{DeePC-Hybrid} from \cite{dorfler2022bridging} and \DeePCsvd{}{} from \cite{zhang2023dimension}. Both of them use two different regularization terms to relax the rank constraint \cref{eqn:bi-level-d} and the row space constraint \cref{eqn:bi-level-e} while \method{DeePC-Hybrid} keeps the Hankel constraint \cref{eqn:bi-level-f} and \DeePCsvd{} drops it. 

Compared with the basic \DeePC{} in \cref{eqn:DeePC}, besides the regularizer $\|\sigma_y\|_2^2$, we introduce two extra regularizers $\|g\|_1$ and $\|(I-\Pi_1)g\|$ in \method{DeePC-Hybrid}, which reads as 
\begin{equation}
\label{eqn:DeePC-Hybrid}
\begin{aligned}
\min_{g, \sigma_y, u \in \mathcal{U}, y \in \mathcal{Y}}  \quad & \|u\|_R^2 + \|y\|_Q^2 + \lambda_1 \|g\|_1 + \lambda_2 \|(I-\Pi_1)g\|_2^2 + \lambda_y \|\sigma_y\|_2^2\\
\textrm{subject to} \quad & 
\begin{bmatrix}
U_\p\\
Y_\p\\
U_\f\\
Y_\f\\
\end{bmatrix} g
= 
\begin{bmatrix}
u_{\ini}\\
y_{\ini}+\sigma_y\\
u\\
y\\
\end{bmatrix} \\
\end{aligned}
\end{equation}
where $\Pi_1 = H_1^\dag H_1$ with $H_1 = \col(U_\p, Y_\p, U_\f)$. Throughout the rest of the discussion, we denote $\|u\|_R^2 = \sum_{k=t}^{t+N-1}\|u(k)\|_R^2$ (similarly for $\|y\|_Q^2$). 
We note that the $l_1$ regularization $\|g\|_1$ can be viewed as a convex relaxation of the rank constraint \cref{eqn:bi-level-e} (\citealp[Theo. IV.8]{dorfler2022bridging}), while the regularization $\|(I-\Pi_1)g\|_2^2$ relaxes row space constraint \cref{eqn:bi-level-d} (\citealp[Theo.~IV.6]{dorfler2022bridging}). 

Since the column number of $H$ is usually larger than its row number in practice (i.e., $H$ is typically a fat matrix), \DeePCsvd{}~in \cite{zhang2023dimension} utilizes singular value decomposition (SVD) to pre-process $H$ and reduce its column dimension. Denoting the SVD of $H$ as $H = W \Sigma V^\tr$, where $\Sigma$ contains its non-zero singular values, we construct a new data matrix $\bar{H} = W \Sigma$, and partition its rows as $\bar{H} = \col(\bar{U}_\p, \bar{Y}_\p, \bar{U}_\f, \bar{Y}_\f)$. Then, the formulation \DeePCsvd{} reads as  
\begin{equation}
\label{eqn:DeePC-SVD}
\begin{aligned}
\min_{\bar{g}, \sigma_y, u \in \mathcal{U}, y \in \mathcal{Y}}  \quad & \|u\|_R^2 + \|y\|_Q^2 + \lambda_1 \|\bar{g}\|_1 + \lambda_2 \|(I-\bar{\Pi}_1)\bar{g}\|_2^2 + \lambda_y \|\sigma_y\|_2^2\\
\textrm{subject to} \quad & 
\begin{bmatrix}
\bar{U}_\p\\
\bar{Y}_\p\\
\bar{U}_\f\\
\bar{Y}_\f\\
\end{bmatrix} \bar{g}
= 
\begin{bmatrix}
u_{\ini}\\
y_{\ini}+\sigma_y\\
u\\
y\\
\end{bmatrix} \\
\end{aligned}
\end{equation}
where $\bar{\Pi}_1 = \bar{H}_1^\dag \bar{H}_1$ and $\bar{H}_1 = \col(\bar{U}_\p, \bar{Y}_\p, \bar{U}_\f)$. 
The dimension of $\bar{g}$ in \cref{eqn:DeePC-SVD} can be much smaller than that in \cref{eqn:DeePC-Hybrid}, and this simple fact can improve numerical efficiency. 

\subsection{Data-Driven Subspace Predictive Control (SPC)}
\label{subsec:data-driven-spc}
We here introduce a new \method{Data-Driven-SPC} to approximate \cref{eqn:bi-level} and establish its equivalence with the classical SPC in \cite{favoreel1999spc}. 
Similar to \method{DeePC-Hybrid} \cref{eqn:DeePC-Hybrid}, we drop the Hankel structure constraint \cref{eqn:bi-level-f} and use a $l_1$ regularization to relax the rank constraint \cref{eqn:bi-level-e}. However, we will directly handle the row space constraint  \cref{eqn:bi-level-d} without using any relaxation. Let us consider
\[
\begin{aligned}
\min_{\tilde{H}} \quad & \|\col(\tilde{U}_\p, \tilde{Y}_\p, \tilde{U}_\f)-\col(U_\p, Y_\p, U_\f)\|  \\
\textrm{subject~to} \quad & \tilde{Y}_\f = Y_\f /\col(\tilde{U}_\p, \tilde{Y}_\p, \tilde{U}_\f).
\end{aligned}
\]
This inner problem has an analytical solution as $\tilde{H}^* = \col(U_\p, Y_\p, U_\f, M) $, with $M = Y_\f \Pi_1$ and $\Pi_1$ defined in \Cref{subsec:DeePC-reg}. Then, we formulate \method{Data-Driven-SPC} as a problem in the form of \cref{eqn:data-driven-spc}. For self-completeness, we also present the classical SPC which is in the form of \cref{eqn:SPC}.   
\vspace{-3pt}

\noindent \begin{minipage}{.48\textwidth} 
\begin{align}
\min_{\sigma_y ,g, u \in \mathcal{U}, y \in \mathcal{Y}}  \;\, & \|u\|_R^2 + \|y\|_Q^2 + \lambda_1 \|g\|_1 + \lambda_y \|\sigma_y\|_2^2 \nonumber \\
\textrm{subject to} \;\, & 
\begin{bmatrix}
U_\p\\
Y_\p\\
U_\f\\
M\\
\end{bmatrix} g
= 
\begin{bmatrix}
u_{\ini}\\
y_{\ini}+\sigma_y\\
u\\
y\\
\end{bmatrix}. \label{eqn:data-driven-spc} 
\end{align}
\end{minipage}
\hspace{1mm}
\begin{minipage}{.48\textwidth}
\vspace{-1mm}
\begin{align}
\min_{\sigma_y, u \in \mathcal{U}, y \in \mathcal{Y}} \;\, & \|u\|_R^2 + \|y\|_Q^2 + \lambda_y \|\sigma_y\|_2^2 \nonumber \\
\textrm{subject to} \;\, & 
y
= Y_\f 
\begin{bmatrix}
    U_\p \\
    Y_p \\
    U_\f
\end{bmatrix}^\dag
\begin{bmatrix}
u_{\ini}\\
y_{\ini}+\sigma_y\\
u\\
\end{bmatrix}. \label{eqn:SPC}
\end{align}
\end{minipage}

\vspace{3pt}

We show that \Cref{eqn:data-driven-spc} is indeed a direct data-driven version of \cref{eqn:SPC} in the sense that they produce the same solution under a very mild condition. The proof is postponed to \Cref{appendix: data-driven-SPC}.  

\begin{thm} \label{proposition:data-driven-SPC}
    If $Q \succ 0, R \succ 0, \lambda_1 = 0$ and $H_1= \col(U_\p, Y_\p, U_\f)$ has full row rank, then \cref{eqn:data-driven-spc} and \cref{eqn:SPC} have the same optimal solution $u^*,y^*,\sigma_y^*$, $\forall \lambda_y >0$.
\end{thm}

Note that our new \method{Data-Driven-SPC} \cref{eqn:data-driven-spc} is more flexible than the classical SPC \cref{eqn:SPC} thanks to the parameter $\lambda_1$, which was motivated from the relaxation of the rank constraint \cref{eqn:bi-level-e}. 

\subsection{\DeePC{} with dominant range space and Hankel structure}
All the convex approximations \cref{eqn:DeePC-Hybrid}, \cref{eqn:DeePC-SVD} and \cref{eqn:data-driven-spc} use different regularizations to relax the difficult constraints \cref{eqn:bi-level-d} and \cref{eqn:bi-level-e}, but both \cref{eqn:DeePC-SVD} and \cref{eqn:data-driven-spc} directly drop the Hankel constraint \cref{eqn:bi-level-f}. In this subsection, we derive another new convex approximation for \cref{eqn:bi-level} that also approximates the Hankel structure with a dominant range space from the SVD. We call it as \method{DeePC-SVD-Iter}. 

\begin{wrapfigure}[5]{r}{.46\textwidth}
\vspace{-5mm}
\begin{subequations} \label{eq:SDV-Dom}
    \begin{align} 
\min_{\tilde{H}} \quad & \|\tilde{H} - H\|_2  \\
\textrm{subject~to} \quad & \textrm{rank}(\tilde{H}) = mL+n  \label{eq:SDV-Dom-a}\\
& \tilde{H} \in \mathcal{H} \label{eq:SDV-Dom-b}
\end{align}
\end{subequations}

\end{wrapfigure}

In particular, we consider \cref{eq:SDV-Dom} as the inner problem in the bi-level formulation \cref{eqn:bi-level}, where the row space constraint \cref{eqn:bi-level-d} will be relaxed using regularization. Note that \cref{eq:SDV-Dom} is also difficult to solve due to the interplay between \cref{eq:SDV-Dom-a} and \cref{eq:SDV-Dom-b}. There are extensive results in the field of structured low-rank approximation (SLRA) \citep{markovsky2008structured}. 
One idea is to use an alternative optimization strategy by considering \cref{eq:SDV-Dom-a} and \cref{eq:SDV-Dom-b} sequentially. Specifically, we here adapt an iterative SLRA algorithm in \cite{yin2021low} to get an approximation solution to \cref{eq:SDV-Dom}. We first note that problem \cref{eq:SDV-Dom} without \cref{eq:SDV-Dom-b} admits an analytical solution $ \tilde{H}^* = W_\R \Sigma_\R V_\R^\tr $
where $W_\R, \Sigma_\R$ and $V_\R$ represent the leading $mL+n$ singular vectors and singular values of $H$, i.e., the dominant range space. 

\begin{wrapfigure}[9]{r}{0.4\textwidth}
\vspace{-6mm}
  \begin{algorithm2e}[H]
    \caption{Iterative SLRA}
    \label{alg:iter-SLRA}
    \KwIn{$H_y$, $\Pi_2$, n, $\epsilon$}
    $H_{y_1} \gets H_y$ 
    
    \Repeat{$\|H_{y_1} -H_{y_2}\| \le \epsilon \|H_{y_1}\|$}{
    $H_{y_2} \gets \hat{H}(H_{y_1})$  \,{\small \texttt{SVD step}}\\
    $H_{y_1} \!\gets\!\! \Pi_\h (H_{y_2})$   {\small \texttt{Hankel~proj}}
    }
    \KwOut{$H_y^* = H_{y_2}$}
\end{algorithm2e}
\end{wrapfigure}

The key idea of the iterative SLRA is to utilize SVD for low-rank approximation of noisy data and then project the low-rank matrix to the set of Hankel matrices. This process is summarized in \Cref{alg:iter-SLRA}. For notational simplicity, We define $H_u = \mathcal{H}_L(u_\D)$, $H_y = \mathcal{H}_L(y_\D)$ to denote the Hankel matrices in \cref{eqn:Hankel_partition}. Thanks to the persistent excitation on $u_\D$ with no noise, we have $\textrm{rank}(H_u) = mL$. However, the measurement $y_\D$ usually contain ``variance'' noise and ``bias'' error, thus the data matrix satisfies $\textrm{rank}(\col(H_u, H_y)) > mL+n$. We use an iterative procedure to denoise $H_y$ while maintaining its Hankel structure (note that we do not change $H_u$ since it is normally noisy-free).  

Let $\Pi_2 = H_u^\dag H_u$ be the orthogonal projector onto the row space of $H_u$, and we first compute $H_y(I -\Pi_2)$ which is the component of $H_y$ in the null space of $H_u$. In each iteration of \Cref{alg:iter-SLRA}, we perform an SVD of  $H_y(I-\Pi_2)$, estimate its rank-$n$ approximation (since we have $\textrm{rank}(H_u) = mL$), and finally combine it with the component of $H_y$ in the row space of $H_u$ as follows 
\[
H_y (I - \Pi_2) = \sum_{i = 1}^{pL} \sigma_i u_i v_i^\tr, \quad \hat{H}(H_y) := H_y \Pi_2 + \sum_{i=1}^{n} \sigma_i u_i v_i^\tr .
\]
We then project $\hat{H}(H_y)$ onto the set of Hankel matrices by averaging skew-diagonal elements and denote  $\Pi_\h$ as the corresponding operator. 
The resulting matrix $H_y^*$ from \Cref{alg:iter-SLRA} is partitioned as $\col(
    Y_\p^*, 
    Y_\f^*) 
$, and we form a new Hankel matrix $\tilde{H}^* = \col(U_\p, Y_\p^*, U_\f, Y_\f^*)$. Finally, we perform an SVD of $\tilde{H}^* =  \tilde{W}_\R \tilde{\Sigma}_\R \tilde{V}_\R^\tr$ to reduce its column dimension and set $\tilde{W}_\R \tilde{\Sigma}_\R = \col(\hat{U}_\p, \hat{Y}_\p, \hat{U}_\f, \hat{Y}_\f)$ with rank $mL + n$ (in the final predictor, we only have the dominant $mL + n$ singular values).   
This new matrix $\col(\hat{U}_\p, \hat{Y}_\p, \hat{U}_\f, \hat{Y}_\f)$ is used as the predictor in \DeePC{} as 
\begin{equation}
\label{eqn:DeePC-SVD-Iter}
\begin{aligned}
\min_{\hat{g}, \sigma_y, u \in \mathcal{U}, y \in \mathcal{Y}}  \quad & \|u\|_R^2 + \|y\|_Q^2  + \lambda_2 \|(I-\hat{\Pi}_1)\hat{g}\|_2^2 + \lambda_y \|\sigma_y\|_2^2\\
\textrm{subject to} \quad & 
\begin{bmatrix}
\hat{U}_\p\\
\hat{Y}_\p\\
\hat{U}_\f\\
\hat{Y}_\f\\
\end{bmatrix} \hat{g}
= 
\begin{bmatrix}
u_{\ini}\\
y_{\ini}+\sigma_y\\
u\\
y\\
\end{bmatrix} \\
\end{aligned}
\end{equation}
where $\hat{\Pi}_1 = \hat{H}_1^\dag \hat{H}_1$, $\hat{H}_1 = \col(\hat{U}_\p, \hat{Y}_\p, \hat{U}_\f)$ and $\|(I-\hat{\Pi}_1)\hat{g}\|_2^2$ is the relaxation term derived from the row space constraint. We call this formulation \cref{eqn:DeePC-SVD-Iter} as \method{DeePC-SVD-Iter}.

\section{Relationship among Different Convex Approximations}
\label{sec:relationships}

As motivated above, \cref{eqn:DeePC-Hybrid,eqn:DeePC-SVD,eqn:DeePC-SVD-Iter,eqn:data-driven-spc} are all tractable convex approximations for the bi-level formulation \cref{eqn:bi-level}. They all begin with the same data matrix $H$ and apply different relaxation strategies to deal with the identification constraints \cref{eqn:bi-level-d,eqn:bi-level-e,eqn:bi-level-f}. We here further look into their relationships and establish certain equivalence. 

First, it is not difficult to see that all of them are equivalent when the data matrix $H$ comes from an LTI system with no noise. We summarize this simple fact below.

\vspace{-2mm}

\begin{fact} \label{fact:1}
    Suppose that the data matrix $H$ in \cref{eqn:Hankel_partition} comes from a controllable LTI system \cref{eqn:LTI-sys} with no noise, and the input $u_\textnormal{d}$ is persistently exciting of order $L+n$. Let $\lambda_1 = 0, \lambda_2 = 0$ and $\sigma_y = 0$. Then, all the \DeePC{} variants in  \cref{eqn:DeePC-Hybrid}, \cref{eqn:DeePC-SVD}, \cref{eqn:data-driven-spc}, and \cref{eqn:DeePC-SVD-Iter} have the same unique optimal solution~$u^*, y^*$. 
\end{fact}

\vspace{-2mm}

For a controllable system \cref{eqn:LTI-sys} with no noise, the data matrix $H$ has already satisfied row space \cref{eqn:bi-level-d}, rank number \cref{eqn:bi-level-e}, and Hankel structure constraints \cref{eqn:bi-level-f}.  Then the new matrix $\tilde{H}^*$ after pre-processing in \method{Data-Driven-SPC} \cref{eqn:data-driven-spc} and \method{DeePC-SVD-Iter} \cref{eqn:DeePC-SVD-Iter} will remain the same as that in \method{DeePC-Hybrid} \cref{eqn:DeePC-Hybrid}, and their range space are also equal to that of \DeePCsvd{}~\cref{eqn:DeePC-SVD}. Thus, all these formulations have the same feasible region and cost functions, and they are equivalent. 

We next move away from noise-free LTI systems. The data matrix $H$ may have `variance'' noise and/or ``bias'' errors; see \cite{dorfler2022bridging}. In this case, we can still show that \method{DeePC-Hybrid} \cref{eqn:DeePC-Hybrid} and \DeePCsvd{} \cref{eqn:DeePC-SVD} produce the same optimal solution $u^*, y^*, \sigma_y^*$.  
\begin{thm}
\label{the:equival-DeePC-SVD}
   Fix any data matrix $H$, and suppose $\lambda_1 = 0$, $\mathcal{U}$ and $\mathcal{Y}$ are convex. Then, \textnormal{\method{DeePC -Hybrid}} \cref{eqn:DeePC-Hybrid} and \textnormal{\DeePCsvd{}} \cref{eqn:DeePC-SVD} have the same optimal solution $u^*, y^*, \sigma_y^*$,  $\forall \lambda_2 >0, \lambda_y>0$. 
\end{thm}

We establish \Cref{the:equival-DeePC-SVD} by expressing $g$ and $\bar{g}$ in terms of $u, y, \sigma_y$. Then, using the SVD properties, we show that \cref{eqn:DeePC-Hybrid} and \cref{eqn:DeePC-SVD} become strictly convex optimization problems with the same objective function, decision variables, and feasible region. The details are presented in  \Cref{appendix:equivalence-DeePC-SVD}. Note that \Cref{the:equival-DeePC-SVD} includes \citealp[Theorem 1]{zhang2023dimension} as a special case, where \cite{zhang2023dimension} requires $\mathcal{U}$ and $\mathcal{Y}$ are convex polytopes that allow simple KKT conditions in their proof.

Finally, \method{DeePC-Hybrid}, \DeePCsvd{}{} and \method{Data-Driven-SPC} are also equivalent under certain conditions. This is summarized in \Cref{prop:spc-equival} below, whose proof is shown in \Cref{appendix: prop:spc-equival}. 
\begin{thm}
\label{prop:spc-equival}
    Fix any data matrix $H$, and suppose $\lambda_1 = 0, \lambda_y >0$ and $\mathcal{U}$ and $\mathcal{Y}$ are convex sets. If $\lambda_2$ is sufficiently large, \textnormal{\method{DeePC-Hybrid}} \cref{eqn:DeePC-Hybrid} and \textnormal{\DeePCsvd{}} \cref{eqn:DeePC-SVD} have the same unique optimal solution $u^*, y^*$ and $\sigma_y^*$ as \textnormal{\method{Data-Driven-SPC}} \cref{eqn:data-driven-spc}.
\end{thm}

The key idea in the proof is to transform the regularizer $\lambda_2 \|I-\Pi_1\|_2^2$ in  \cref{eqn:DeePC-Hybrid} as the constraint $\|I-\Pi_1\|_2 = 0$ when $\lambda_2$ is sufficiently large via penalty arguments. Then, \cref{eqn:DeePC-Hybrid} and \cref{eqn:data-driven-spc} have the same objective function and decision variables. The proof is completed by further establishing that they have the same feasible region. From \Cref{proposition:data-driven-SPC,the:equival-DeePC-SVD,prop:spc-equival}, we conclude that \method{Data-Driven-SPC} \cref{eqn:data-driven-spc}, \method{DeePC-Hybrid} \cref{eqn:DeePC-Hybrid}, \DeePCsvd{}{} \cref{eqn:DeePC-SVD} are equivalent to classical SPC \cref{eqn:SPC} with noisy data when $\lambda_1 = 0, \lambda_y >0$, $\lambda_2$ is sufficiently large and $H_1$ has full row rank. This is more general~than \cite{fiedler2021relationship,breschi2023data}: the equivalence for \method{DeePC-Hybrid} with regularizer $\|g\|_2^2$ and classical SPC is discussed in \cite{fiedler2021relationship}, while \method{DeePC-Hybrid} and an approach similar to \method{Data-Driven-SPC} are proved to be equivalent in \cite{breschi2023data}.       

Note that the new variant \method{DeePC-SVD-Iter} involves an iterative algorithm to pre-process the noisy data (\Cref{alg:iter-SLRA}), and thus it is non-trivial to formally establish its relationship with respect to other variants. Yet, our numerical experiments in  \cref{sec:results} show that \method{DeePC-SVD-Iter} often has superior performance among all these convex approximations for noisy data.  

\section{Numerical Experiments}
\label{sec:results}
We perform numerical experiments to illustrate \Cref{the:equival-DeePC-SVD,prop:spc-equival}\footnote{Our code is available at \url{https://github.com/soc-ucsd/Convex-Approximation-for-DeePC}.}.  We also numerically investigate the effects of $\lambda_1, \lambda_2$, and confirm the superior performance of \method{DeePC-SDV-iter} \cref{eqn:DeePC-SVD-Iter}. Some additional numerical results on nonlinear systems are provided in \Cref{appendix:numerical-results}.

\begin{figure}[t]
\setlength{\abovedisplayskip}{0pt}
\centering
\subfigure[]{\includegraphics[width=0.32\textwidth]{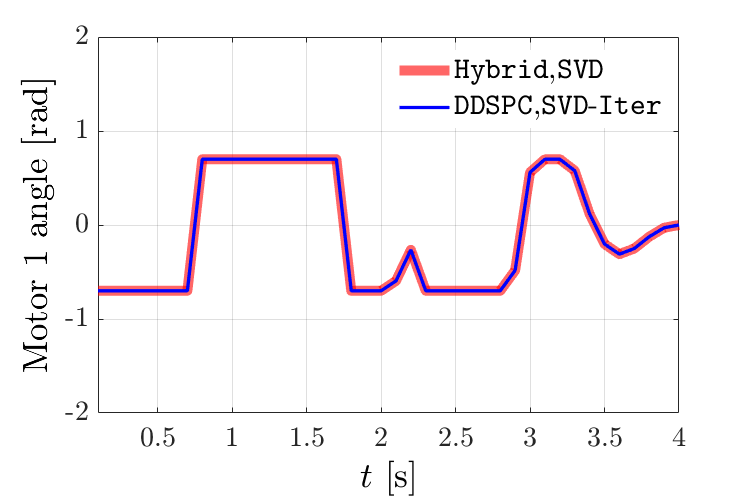} \label{fig:eq-all}}
\subfigure[]{\includegraphics[width=0.32\textwidth]{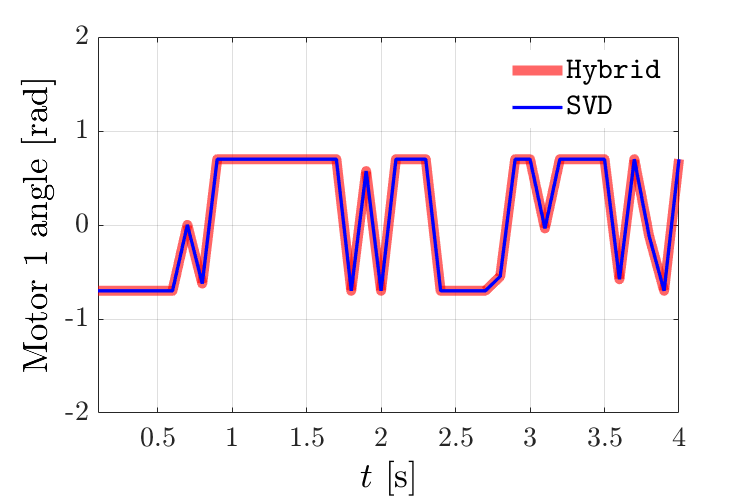}\label{fig:eq-hy-svd}}
\subfigure[]{\includegraphics[width=0.32\textwidth]{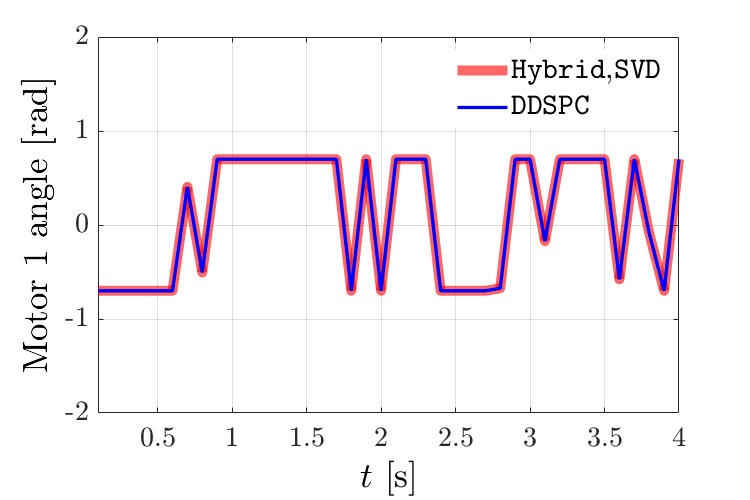}\label{fig:eq-hy-svd-spc}}

\vspace{-4mm}

\caption{Equivalent optimal solutions of different convex approximations in \Cref{the:equival-DeePC-SVD,prop:spc-equival}. (a) All methods with noise-free data. (b) \method{DeePC-Hybrid} and \DeePCsvd{} with $\lambda_1\! =\! 0, \lambda_2\!=\!30$ and $\lambda_y\! =\! 100$. (c) \method{DeePC-Hybrid}, \DeePCsvd{} and \method{Data-Driven-SPC} with $\lambda_1\! =\! 0, \lambda_2\!=\!10000$ and $\lambda_y\! =\! 100$.
}
\vspace{-4mm}
\end{figure}

\noindent \textbf{Experiment setup.}  We consider an LTI system from \cite{fiedler2021relationship}. This is a triple-mass-spring system with $n=8$ states, $m=2$ inputs (two stepper motors), and $p=3$ outputs (disc angles). In our experiments, the length of the pre-collected trajectory is $T=200$, and the prediction horizon and the initial sequence are chosen as $N = 40$ and $T_{\ini} = 4$, respectively. We choose $Q = I$, $R = 0.1I$ and  $\mathcal{U} = [-0.7, 0.7]$.

\noindent \textbf{Equivalence.}  We here numerically verify that the optimal solutions from different convex approximations are the same under appropriate settings. For noise-free pre-collected data (\Cref{fact:1}), all methods have the same optimal solution, and one solution instance is given in  \Cref{fig:eq-all}. We next consider data collection with additive Gaussian measurement noises $\omega \sim \mathcal{N}(0, 0.01 I)$. 
We choose $\lambda_1 = 0, \lambda_2 = 30$ and $\lambda_y = 100$ according to \Cref{the:equival-DeePC-SVD}. One solution instance is shown in \Cref{fig:eq-hy-svd}, which shows that \method{DeePC-Hybrid} and \DeePCsvd{} provide the same optimal solution. Finally, for the equivalence of \method{Data-Driven-SPC}, \method{DeePC-Hybrid} and \DeePCsvd{}, we choose $\lambda_1 = 0$ and {$\lambda_2 = 10000$} according to \Cref{prop:spc-equival}, and the results are shown in \Cref{fig:eq-hy-svd-spc}.

\begin{wrapfigure}[13]{r}{.46\textwidth}
\vspace{-5.5mm}
\centering
\includegraphics[width=0.45 \textwidth]{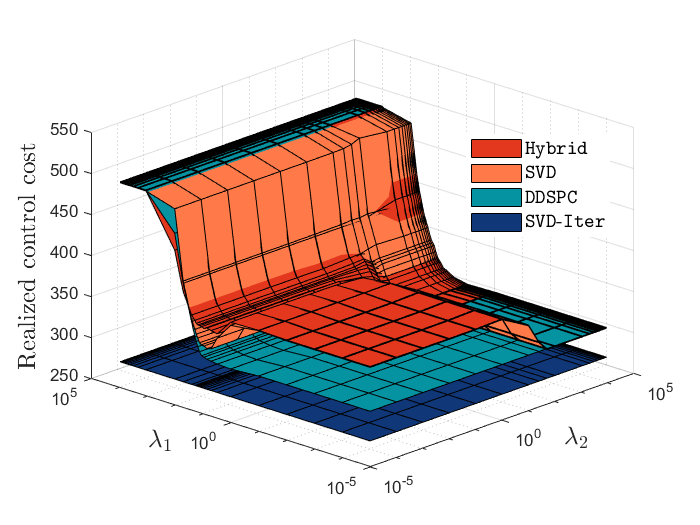}
\vspace{-5mm}
\caption{Realized control cost where sampling points are chosen from $\lambda_1, \lambda_2 \in [10^{-5},10^{4}]$ to capture the effect of hyperparameters.}
\label{fig:hyper-para}
\end{wrapfigure}

\noindent \textbf{Influence of $\lambda_1$ and $\lambda_2$.} We then analyze the effect of hyperparameters $\lambda_1$ and $\lambda_2$. In particular, similar to \cite{dorfler2022bridging}, we consider the realized control cost after applying the optimal control inputs from  \cref{eqn:DeePC-Hybrid,eqn:DeePC-SVD,eqn:data-driven-spc,eqn:DeePC-SVD-Iter}, which is computed as $\|u_\textrm{opt}\|_R^2+\|y_\textrm{true}\|_Q^2$, where $u_\textrm{opt}$ is the computed optimal control input and $y_\textrm{true}$ is the realized trajectory after applying it. We fixed $\lambda_y = 100$.

\Cref{fig:hyper-para} shows the realized control performance over $\lambda_1$ and $\lambda_2$ for different convex approximations. The hyperparameters $\lambda_1$ and $\lambda_2$  indeed have a significant effect for \method{DeePC-Hybrid} \cref{eqn:DeePC-Hybrid}, \DeePCsvd{} \cref{eqn:DeePC-SVD} and \method{Data-Driven-SPC} \cref{eqn:data-driven-spc} (denoted as \method{DDSPC} in \Cref{fig:hyper-para}). For these methods, $\lambda_1$ needs to be chosen more carefully, neither too large nor too small and $\lambda_2$ should not be chosen too small; similar phenomena also appeared in \cite{dorfler2022bridging}. However, it is notable that \method{DeePC-SVD-Iter} \cref{eqn:DeePC-SVD-Iter} not only achieves the best performance but also is not very sensitive to $\lambda_2$ (note that $\lambda_1 = 0$ in \cref{eqn:DeePC-SVD-Iter}).

\noindent \textbf{Comparison of \DeePC{} variants.} Finally, we compare the realized control cost and the computational time for different convex approximations. Motivated by \Cref{fig:hyper-para}, we choose $\lambda_1 = \lambda_2 = 30$ and $\lambda_y = 100$. The performance of \DeePC{} variants is related to the pre-collected trajectory. Thus, all presented realized control costs and computational time for different convex approximations are averaged over 100 pre-collected trajectories. The numerical results are listed in \Cref{table:cost-time}. The ground-truth cost is computed from \cref{eqn:DeePC} with noise-free data. From \Cref{table:cost-time}, we see that the realized control cost satisfies $\method{DeePC-Hybrid}>\DeePCsvd{}>\method{Data-Driven-SPC}>\method{DeePC-SVD-Iter}>\textrm{System ID}$. For the LTI system with noisy data, the inner problem in \Cref{eqn:bi-level} forces the data-driven representation to be more structured, which enhances noise rejection performance for upper predictive control in \Cref{eqn:bi-level}. The increasing rate of realized cost for our new \method{DeePC-SVD-Iter} is $3.9\%$, which is much better than other \DeePC{} variants. 
 
\begin{figure}[t]
\centering
\subfigure[\small \method{Hybrid}]{\includegraphics[width=0.2\textwidth]{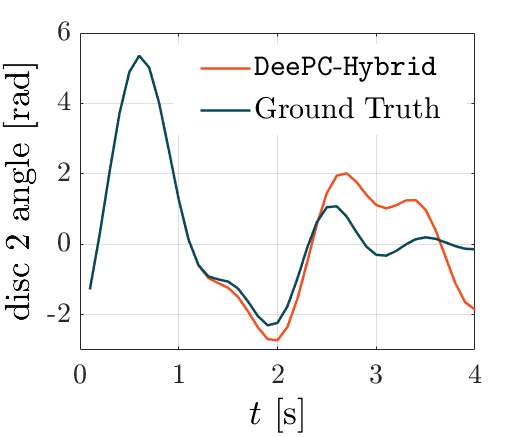} \label{fig:traj-hybrid}}
\hspace{-2mm}\subfigure[\small \method{SVD}]{\includegraphics[width=0.2\textwidth]{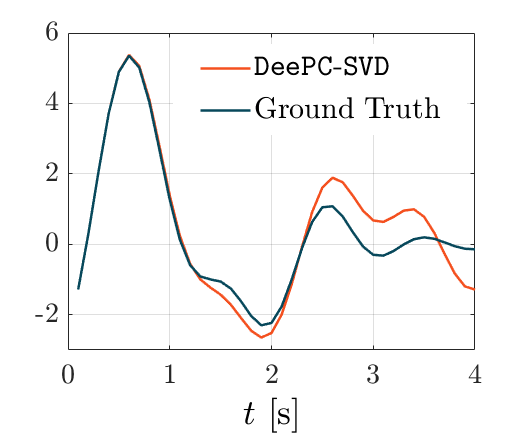}}
\hspace{-2mm}\subfigure[\small\method{DDSPC}]{\includegraphics[width=0.2\textwidth]{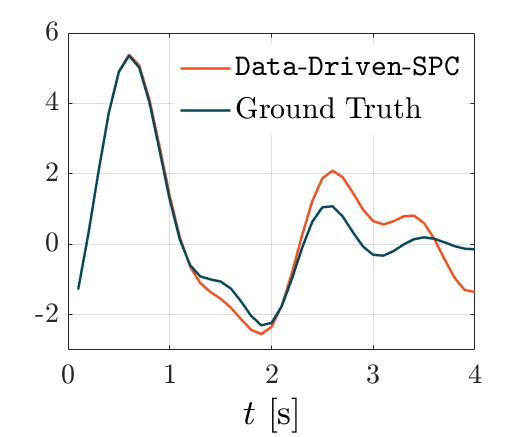}}
\hspace{-2mm}\subfigure[\small \method{SVD-Iter}]{\includegraphics[width=0.2\textwidth]{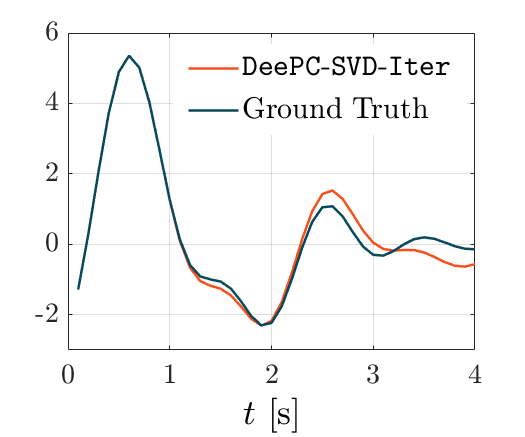}}
\hspace{-2mm}\subfigure[\small System ID]
{\includegraphics[width=0.2\textwidth]{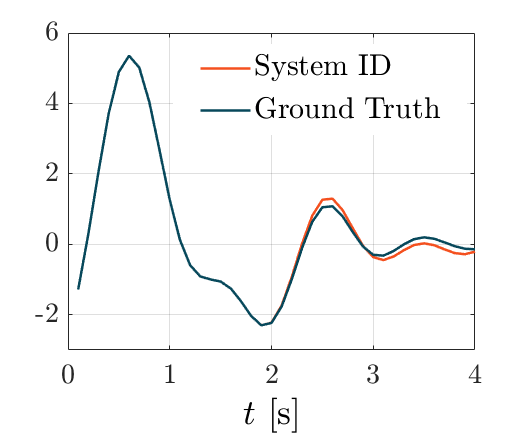} \label{fig:traj-sysID}}
\vspace{-3mm}
\caption{Open-loop trajectories (the angle of disc 2) of different methods. The blue trajectory and orange trajectories represent ground truth and different approximation methods, respectively.
}
\label{fig:traj-nond-LTI}

\vspace{-2.5mm}
\end{figure}

\Cref{fig:traj-nond-LTI} shows one typical open-loop trajectory for all methods. In this case, the open-loop trajectories from \cref{eqn:DeePC-Hybrid,eqn:DeePC-SVD,eqn:data-driven-spc,eqn:DeePC-SVD-Iter} remain close to the ground truth up to 2 $\mathrm{s}$. Then the trajectory is better aligned with the ground truth from Fig. \ref{fig:traj-hybrid} to Fig. \ref{fig:traj-sysID} as the corresponding data-driven representation becomes more structured. Our numerical results also suggest that the indirect system ID approach is superior in the case of ``variance'' noise, consistent with \cite{dorfler2022bridging}. In \Cref{appendix:numerical-results}, our numerical results on nonlinear systems further reveal that \method{DeePC-SVD-Iter} \cref{eqn:DeePC-SVD-Iter} also has enhanced performance in the case of ``bias'' errors. 

\begin{table}[t]
\setlength{\abovedisplayskip}{0pt}
\caption{\vspace{-1mm} Realized Control Cost and Computational Time; GT denotes ground truth with noisy-free data. }
\vspace{-2mm}
\centering
{\footnotesize
\begin{threeparttable}
    \begin{tabular}{ccccccc}
    \toprule
    &GT&\method{Hybrid} & \method{SVD} & \method{Data-Driven-SPC} & \method{SVD-Iter} & System ID\\
    \midrule
    Realized Cost & 277.25 & 388.42 & 370.20 & 365.08  & 288.17 & 279.86  \\
    Increasing Rate &  N/A & $40.1 \%$ & $33.5\%$ & $31.7\%$ & $\mathbf{3.9\%}$ & $0.9 \%$\\
    Compu. time [$\mathrm{s}$] & N/A & 0.133 & 0.131 & 0.097 & 0.104 & N/A \\
    \bottomrule
    \end{tabular}
    \end{threeparttable}
    }
    \label{table:cost-time}
    \vspace{-1mm}
\end{table}

\vspace{-3mm}

\section{Conclusion}
\label{sec:conclusions}
In this paper, we have proposed a new bi-level formulation incorporating system ID techniques and predictive control. The existing \DeePC{} (i.e., \method{DeePC-Hybrid} \cref{eqn:DeePC-Hybrid} and \method{DeePC-SVD} \cref{eqn:DeePC-SVD}) and also new variants (i.e., \method{Data-driven-SPC} \cref{eqn:data-driven-spc} and \method{DeePC-SVD-Iter} \cref{eqn:DeePC-SVD-Iter}) can be considered as convex approximations of this bi-level formulation. We have further clarified their equivalence under appropriated settings (\Cref{the:equival-DeePC-SVD,prop:spc-equival,proposition:data-driven-SPC}). Numerical simulations have validated our theoretical findings, and also revealed the superior performance of \method{DeePC-SVD-Iter} \cref{eqn:DeePC-SVD-Iter} with a more structured predictor. Interesting future directions include analyzing the effect of the length of pre-collected data, and investigating the closed-loop performance of different \DeePC{} variants.


\acks{This work is supported in part by NSF ECCS-2154650 and in part by NSF CMMI-2320697.}

\bibliography{reference.bib}

\newpage
\appendix

\section{Technical proofs}

In this section, we provide the technical proofs for  \Cref{proposition:data-driven-SPC,the:equival-DeePC-SVD,prop:spc-equival} that are omitted in the main text. 

\subsection{Linear algebra fundamentals} \label{appendix:review}

Before presenting our proof details, let us briefly refresh some standard linear algebra, including singular value decomposition (SVD), pseudo inverse, and orthogonal decomposition. Essentially, many topics related to \method{DeePC} are all boiled down to some linear algebra. 

We consider real matrices $A \in \mathbb{R}^{m_\A \times n_\A}$ and $B \in \mathbb{R}^{m_\B \times n_\B}$. The row space of $A$ is denoted as $\textrm{rowsp}(A)$ which is the span of its row vectors while the column space of $A$ is defined as the span of its column vectors. Let rank$(A) = r$, the compact singular value decomposition (SVD) of $A$~is
\begin{equation}
\label{eqn:SVD-compact}
A = W \Sigma V^\tr
\end{equation}
where $\Sigma = \textrm{diag}(\sigma_1, \ldots, \sigma_r)\in \mathbb{R}^{r \times r}$, $W = [w_1,\ldots, w_r] \in \mathbb{R}^{m_\A\times r}$, $V = [v_1, \ldots, v_r] \in \mathbb{R}^{n_\A \times r}$ and $\sigma_i$, $w_i$ and $v_i$ are singular values, left singular vectors and right singular vectors of $A$, respectively. Matrices $\Sigma, W, V$ also satisfies $\sigma_1 \ge \cdots \ge \sigma_r >0$ and $W^\tr W=I$, $V^\tr V=I$ (i.e., $w_i$ and $v_i$ are othornormal vectors). We note that the range space and row space of $A$ are the same as the range space of $W$ and $V$, respectively.

The pseudo inverse of a matrix $A$, represented as $A^\dag$, is a generalization of the inverse matrix. It is commonly used to solve a system of linear equations $Ax = b$.  The pseudo inverse provides a least-squares approximated solution $x = A^\dag b$ if there is no exact solution for the original problem. On the other hand, if the solution exists, $x_p=A^\dag b$ is a least-norm solution and the general solution can be represented as $x = x_p + \hat{x}$ where $\hat{x}$ is a vector in the null space of $A$. The pseudo inverse matrix $A^\dag$ needs to satisfy the following four criteria \citep{penrose1955generalized}
\begin{equation} \label{eq:pseudo-inverse}
A A^\dag A = A, \qquad A^\dag A A^\dag = A^\dag, \qquad (AA^\dag)^\tr = AA^\dag, \qquad (A^\dag A)^\tr = A^\dag A, 
\end{equation}
which is unique and can be computed by compact SVD of $A$ in \cref{eqn:SVD-compact} as $A^\dag = V \Sigma^{-1} W^\tr$. By definition, it is clear that $AA^\dag$ and $A^\dag A$ are symmetric matrices. We list the following properties of pseudo inverse that will be utilized in our proofs later \citep{golub2013matrix,greville1966note}:
\begin{subequations}
\begin{equation}
\label{eqn:p-full-row-rank}
     \text{If $A$ has full row rank, then $A^\dag = A^\tr(AA^\tr)^{-1}$ and $AA^\dag = I$}.
\end{equation}
     \begin{equation}
\label{eqn:p-full-col-rank}
     \text{If $A$ has full column rank, then $A^\dag = (A^\tr A)^{-1}A^\tr$ and $A^\dag A = I$}.
\end{equation}
     \begin{equation}
\label{eqn:p-range}
     \text{The range space of $A^\dag$ is the same as the row space of $A$}.
\end{equation}
     \begin{equation}
\label{eqn:p-orthonormal}
     \text{If column vectors of $A$ are orthonormal, we have $(A^\tr)^\dag = A$}.
\end{equation}
     \begin{equation}
\label{eqn:p-exchange}
     \text{If $B$ has orthonormal rows, we have $(AB)^\dag = B^\dag A^\dag$}.
\end{equation}
\end{subequations}

The orthogonal decomposition of a vector $u \in \mathbb{R}^n$ with respect to a subspace $F \subseteq \mathbb{R}^n$ is 
$
u = v + v^\bot,
$
where $v \in F$ and $v^\bot$ in its orthogonal complement $F^\bot$, i.e., $v^\tr v^\bot = 0$ \citep{golub2013matrix}. The vector $v$ is called the orthogonal projection of $u$ onto $F$. Given a matrix $A=[a_1, \ldots, a_{m_\A}]^\tr$ (i.e., $a_i^\tr$ is a row vector of $A$), its row space orthogonal projection  onto the subspace $F$, denoted as $M = [m_1, \ldots, m_{m_\A}]^\tr$, is constructed by the orthogonal projection of $a_i$ onto $F$ via  $a_i = m_i + m_i^\bot, i=1,\ldots,m_\A$ (the column space orthogonal projection is similar). The orthogonal projector onto the {range space of $A^\tr$ (or equivalently row space of $A$) is $A^\dag A$}, and the orthogonal decomposition of $u=v + v^\bot$ associated with range space of $A^\tr$ can be computed~as
\[
v = A^\dag A u, \quad v^\bot = (I-A^\dag A)u,
\]
and the row space orthogonal projection $M$ of matrix $B$ onto the row space of $A$ is denoted as  
\[
M = B/A = ((A^\dag A)B^\tr)^\tr = B(A^\dag A). 
\]

\subsection{Proof of \Cref{proposition:data-driven-SPC}} \label{appendix: data-driven-SPC}

Our proof is divided into two main parts:
\begin{enumerate}
    \item When $H_1= \col(U_\p, Y_\p, U_\f)$ has full row rank, we show that \cref{eqn:data-driven-spc} and \cref{eqn:SPC} have the same feasible region: if $g, u, y, \sigma_y$ is feasible to  \cref{eqn:data-driven-spc}, then the same $u, y, \sigma_y$ is also feasible for \cref{eqn:SPC}. Conversely, given any feasible solution $u, y, \sigma_y$ to \cref{eqn:SPC}, we can construct a vector $g$ such that $g, u, y, \sigma_y$ is feasible to  \cref{eqn:data-driven-spc}.
    \item When $\lambda_1 = 0$, \cref{eqn:data-driven-spc} and \cref{eqn:SPC} have the same cost function in terms of $u, y, \sigma_y$. 
\end{enumerate}
Combining the two properties above with the fact that the cost function in \cref{eqn:SPC} is strongly convex, we conclude that \cref{eqn:data-driven-spc} and \cref{eqn:SPC} have the same unique optimal solution $u^*, y^*, \sigma_y^*$. 

The property 2 above is obvious. We prove the property 1 below. Let us first decompose 
     \begin{equation} \label{eq:ortho-decomposition-Yf}
     Y_\f = M+M^\bot
     \end{equation}
     where $M$ is the (row space) orthogonal projection of $Y_\f$ on the row space of $H_1$ and $M^\bot$ is the rest part of $Y_\f$ in the null space of $H_1$. Since $H_1$ has full row rank, we have $H_1 H_1^\dag = I$ from \cref{eqn:p-full-row-rank}. Also, the range space of $H_1^\dag$ is the same as the row space of $H_1$ (see \cref{eqn:p-range}), which means $M^\bot H_1^\dag = 0$. 
    
    We assume $u_1, y_1, \sigma_{y_1}, g_1$ is a feasible solution for \cref{eqn:data-driven-spc}. Then, without loss of generality, the vector $g_1$ can be represented as 
    \[
    g_1 = H_1^\dag \col(
        u_\ini, y_\ini + \sigma_{y_1}, u_1
    ) + \hat{g}
    \]
    where $\hat{g}$ is a vector in the null space of $H_1$. We have $M \hat{g} = 0$ since $H_1\hat{g} = 0$ and $Y_\f H_1^\dag = (M+M^\bot)H_1^\dag = M H_1^\dag$ because 
    $M^\bot H_1^\dag = 0$. Thus, from the equality constrain in \cref{eqn:data-driven-spc}, the vector $y_1$ satisfies
    \begin{equation}
    \label{eqn:relate-y-YF}
    y_1 = M g_1 = M H_1^\dag \begin{bmatrix}
        u_\ini \\ y_\ini + \sigma_{y_1} \\ u_1
    \end{bmatrix} + M \hat{g} = M H_1^\dag \begin{bmatrix}
        u_\ini \\ y_\ini + \sigma_{y_1} \\ u_1
    \end{bmatrix} 
    = Y_\f H_1^\dag \begin{bmatrix}
        u_\ini \\ y_\ini + \sigma_{y_1} \\ u_1
    \end{bmatrix},
    \end{equation}
    which means $u_1, y_1, \sigma_{y_1}$ is also a feasible solution of $\cref{eqn:SPC}$.
    
    We next assume $u_1, y_1, \sigma_{y_1}$ is a feasible solution for \cref{eqn:SPC}. Substituting the orthonormal decomposition \cref{eq:ortho-decomposition-Yf} into the equality constraint of \cref{eqn:SPC}, we have
    \begin{equation}
    \label{eqn:Yf_derive}
    y_1 = Y_\f H_1^\dag \begin{bmatrix}
        u_\ini \\ y_\ini + \sigma_{y_1} \\ u_1
    \end{bmatrix} = (M+M^\bot) H_1^\dag \begin{bmatrix}
        u_\ini \\ y_\ini + \sigma_{y_1} \\ u_1
    \end{bmatrix} 
     = M H_1^\dag \begin{bmatrix}
        u_\ini \\ y_\ini + \sigma_{y_1} \\ u_1
    \end{bmatrix}.
    \end{equation}
    Upon defining $g_1 = H_1^\dag \col(
        u_\ini, y_\ini + \sigma_{y_1}, u_1
    )$, we have $y_1 = M g_1$ from \cref{eqn:Yf_derive}. We then substitute $g_1$ into the equality constraint of \cref{eqn:data-driven-spc}, leading to
    \[
    \begin{bmatrix}
        H_1 \\ M
    \end{bmatrix} g_1 = 
    \begin{bmatrix}
    H_1 H_1^\dag \begin{bmatrix}
        u_\ini \\ y_\ini + \sigma_{y_1} \\ u_1
    \end{bmatrix}  \\
    M g_1
    \end{bmatrix}
    = \begin{bmatrix}
        u_\ini \\ y_\ini + \sigma_{y_1} \\ u_1 \\
        y_1
    \end{bmatrix},
    \]
where we have used the fact $H_1 H_1^\dag = I$ from \cref{eqn:p-full-row-rank} since $H_1$ has full row rank. This means that $g_1, u_1, y_1, \sigma_{y_1}$ is a feasible solution for \cref{eqn:data-driven-spc}. This completes our proof.   

\subsection{Proof of \Cref{the:equival-DeePC-SVD}} \label{appendix:equivalence-DeePC-SVD} 
Since $\lambda_1 = 0$, we only consider a two-norm regularizer $\|\cdot\|_2$. Here, we consider \method{DeePC-Hybrid} \cref{eqn:DeePC-Hybrid} and \DeePCsvd{}{} \cref{eqn:DeePC-SVD} with a general two-norm regularization. For convenience, we rewrite their forms below:
\begin{equation}
\label{eqn:DeePC-Hybrid-2-norm}
\begin{aligned}
\min_{g, \sigma_y, u \in \mathcal{U}, y \in \mathcal{Y}}  \quad & \|u\|_R^2 + \|y\|_Q^2 + \lambda_2 \|Gg\|_2^2 + \lambda_y \|\sigma_y\|_2^2\\
\textrm{subject to} \quad & 
\begin{bmatrix}
U_\p\\
Y_\p\\
U_\f\\
Y_\f\\
\end{bmatrix} g
= 
\begin{bmatrix}
u_{\ini}\\
y_{\ini}+\sigma_y\\
u\\
y\\
\end{bmatrix},
\end{aligned}
\end{equation}
\begin{equation}
\label{eqn:DeePC-SVD-2-norm}
\begin{aligned}
\min_{\bar{g}, \sigma_y, u \in \mathcal{U}, y \in \mathcal{Y}}  \quad & \|u\|_R^2 + \|y\|_Q^2 + \lambda_2 \|\bar{G}\bar{g}\|_2^2 + \lambda_y \|\sigma_y\|_2^2 \\
\textrm{subject to} \quad & 
\begin{bmatrix}
    \bar{U}_\p \\
    \bar{Y}_\p \\
    \bar{U}_\f \\
    \bar{Y}_\f
\end{bmatrix}
\bar{g}
= 
\begin{bmatrix}
u_{\ini}\\
y_{\ini}+\sigma_y\\
u\\
y\\
\end{bmatrix}. \\
\end{aligned}
\end{equation}
We recall that the SVD decomposition: 
$$
\begin{aligned}
H &= \col(U_\p, Y_\p, U_\f, Y_\f) = W \Sigma V^\tr, \\
\bar{H} &= \col(\bar{U}_\p, \bar{Y}_\p, \bar{U}_\f, \bar{Y}_\f) = W \Sigma, 
\end{aligned}$$
where $\Sigma$ contains $r$ nonzero singular values, the columns of $W\Sigma$ are linearly independent and the columns of $V$ are orthonormal vectors. The orthonormal columns lead to that $V^\tr V = I$ and $(V^\tr)^\dag = V$ (see \cref{eqn:p-full-col-rank} and \cref{eqn:p-orthonormal}). We note that thanks to the SVD, the variable $\bar{g}$ has a dimension of $r$ which is typically much smaller than the dimension of $g$. 

We next provide a set of sufficient conditions for the equivalence of \cref{eqn:DeePC-Hybrid-2-norm} and \cref{eqn:DeePC-SVD-2-norm}. Its proof is provided in \Cref{appendix:proposition}. 
\begin{prop}
\label{prop:reg-equil-condtion}
Let $Q\succ 0, R \succ 0$  in \cref{eqn:DeePC-Hybrid-2-norm} and \cref{eqn:DeePC-SVD-2-norm}. If the matrices $G \in \mathbb{R}^{m_g \times (T-T_\ini-N+1)}, \bar{G}\in \mathbb{R}^{m_{\bar{g}} \times r}, 
V\in \mathbb{R}^{(T-T_\ini-N+1) \times r}$ and $H\in \mathbb{R}^{(m+p)(T_\ini+N)\times (T-T_\ini-N+1)}$ satisfy two properties below
\begin{subequations}
    \begin{align}
    \textnormal{rowsp}(H G^\tr G) &\subseteq  \textnormal{rowsp}(H), \label{eqn:prop-1} \\
    V^\tr G^\tr G V & = \bar{G}^\tr \bar{G} \label{eqn:prop-2}, 
    \end{align}
\end{subequations}
then the optimal solution for $u, y, \sigma_y$ of \cref{eqn:DeePC-Hybrid-2-norm} and \cref{eqn:DeePC-SVD-2-norm} are the same and unique. 
\end{prop}

Then, \Cref{the:equival-DeePC-SVD} can be considered as a corollary of \Cref{prop:reg-equil-condtion}. 

\noindent \textbf{Proof of \Cref{the:equival-DeePC-SVD}:} 
It suffices to show that $G = I -\Pi_1$ and $\bar{G} = I- \bar{\Pi}_1$ satisfy the two properties \cref{eqn:prop-1} and \cref{eqn:prop-2} in \Cref{prop:reg-equil-condtion}. 

For simplicity, let us denote $v=\col(u_\ini,y_\ini+\sigma_y,u,y)$. 
We also define $\bar{H}_1 = \col(\bar{U}_\p, \bar{Y}_\p, \bar{U}_\f)$ and recall that $H_1 = \col(U_\p, Y_\p, U_\f)=\bar{H}_1 V^\tr$ due to the SVD decomposition. Then, the matrices $G$ and $\bar{G}$ can be represented as 
\[
\begin{aligned}
    G &= I - H_1^\dag H_1 = I - (\bar{H}_1 V^\tr)^\dag \bar{H}_1 V^\tr = I - V \bar{H}_1^\dag \bar{H}_1 V^\tr,  \\ 
    \bar{G} &= I -\bar{H}_1^\dag \bar{H}_1, 
\end{aligned}
\]
in which we have used the fact $(\bar{H}_1 V^\tr)^\dag \! =\! V \bar{H}_1^\dag$ from \cref{eqn:p-exchange} since $V^\tr$ has orthonormal row~vectors. 

To establish \cref{eqn:prop-1}, we first have the following relationship 
\begin{equation} \label{eq:HGG}
\begin{aligned}
H G^\tr G &= H (I-H_1^\dag H_1)^\tr(I-H_1^\dag H_1) \\
&= H - H H_1^\tr \times (H_1^\dag)^\tr - (H - H H_1^\tr (H_1^\dag)^\tr)H_1^\dag \times H_1.
\end{aligned}
\end{equation}
Note that for any pair of compatible matrices, we have 
\begin{equation} \label{eq:AB-B}
\textrm{rowsp}(AB) \subseteq \textrm{rowsp}(B).    
\end{equation}
Combining \cref{eq:HGG} with \cref{eq:AB-B}, to establish $\textnormal{rowsp}(H G^\tr G) \subseteq  \textnormal{rowsp}(H)$, it suffices to prove that the row spaces of $(H_1^\dag)^\tr$ and $H_1$ are subspaces of the row space of $H$. Since $H_1$ is constructed by the top three block matrices of $H$, we naturally have $\textrm{rowsp}(H_1) \subseteq \textrm{rowsp}(H)$. Meanwhile, the property \cref{eqn:p-range} of the pseudo inverse directly implies that $\textrm{rowsp}((H_1^\dag)^\tr)=\textrm{rowsp}(H_1) \subseteq \textrm{rowsp}(H)$.

For \cref{eqn:prop-2}, we have
\[
\begin{aligned}
V^\tr G^\tr G V &= V^\tr (I-V \bar{H}_1^\tr (\bar{H}_1^\dag)^\tr V^\tr)(I - V \bar{H}_1^\dag \bar{H}_1 V^\tr)V \\
&= (I- \bar{H}_1^\tr (\bar{H}_1^\dag)^\tr )V^\tr V(I - \bar{H}_1^\dag \bar{H}_1 ) = \bar{G}^\tr \bar{G},
\end{aligned}
\]
where we have used the fact that $V^\tr V = I$. 
This means \cref{eqn:prop-2} is also satisfied. 

Then, \Cref{the:equival-DeePC-SVD} becomes a corollary of \Cref{prop:reg-equil-condtion}. This completes our proof.

\subsection{Proof of \Cref{prop:spc-equival}} \label{appendix: prop:spc-equival}
We present the proof of \Cref{prop:spc-equival} for \method{DeePC-Hybrid} since \method{DeePC-Hybrid} and \DeePCsvd{} are equivalent when $\lambda_1 = 0$ (see \Cref{the:equival-DeePC-SVD}). The key idea is that we can  write \method{DeePC-Hybrid} equivalently as
\begin{subequations}
\label{eqn:equival-data-driven-SPC}
\begin{align}
\min_{g, \sigma_y, u \in \mathcal{U}, y \in \mathcal{Y}}  \quad & \|u\|_R^2 + \|y\|_Q^2 + \lambda_y \|\sigma_y\|_2^2 \nonumber \\
\textrm{subject to} \quad & 
\begin{bmatrix}
U_\p\\
Y_\p\\
U_\f\\
Y_\f\\
\end{bmatrix} g
= 
\begin{bmatrix}
u_{\ini}\\
y_{\ini}+\sigma_y\\
u\\
y\\
\end{bmatrix}, \label{eqn:equival-data-driven-SPC-1}\\
& \|(I-\Pi_1)g\|_2 = 0, \label{eqn:equival-data-driven-SPC-2}
\end{align}
\end{subequations}
when $\lambda_2$ is sufficiently large \cite[Proposition. IV.2 \& Proposition. A.3]{dorfler2022bridging}. It is obvious that \cref{eqn:equival-data-driven-SPC} and \cref{eqn:data-driven-spc} have the same objective function when $\lambda_1=0$ and it only contains $u, y$ and $\sigma_y$. Thus, we show that \cref{eqn:equival-data-driven-SPC} and \cref{eqn:data-driven-spc} provide the same unique optimal solution $u^*, y^*$ and $\sigma_y^*$ by proving:
\begin{enumerate}
    \item Feasible regions of $u,y$ and $\sigma_y$ are the same for \cref{eqn:equival-data-driven-SPC} and \cref{eqn:data-driven-spc}: if $g, u, y, \sigma_y$ is feasible to  \cref{eqn:equival-data-driven-SPC}, then the same $g, u, y, \sigma_y$ is also feasible for \cref{eqn:data-driven-spc}. Conversely, given any feasible solution $g, u, y, \sigma_y$ to \cref{eqn:data-driven-spc}, we can construct a vector $\tilde{g}$ such that $\tilde{g}, u, y, \sigma_y$ is feasible to \cref{eqn:equival-data-driven-SPC}.
    \item The optimal solution of $u, y$ and $\sigma_y$ is unique for \cref{eqn:data-driven-spc}. 
\end{enumerate}

We assume $u_1, y_1, \sigma_{y_1}$ and $g_1$ is a feasible solution for \cref{eqn:equival-data-driven-SPC}. Substituting the orthogonal decomposition \cref{eq:ortho-decomposition-Yf} of $Y_\f$ into \cref{eqn:equival-data-driven-SPC-1}, we have
\[
\begin{bmatrix}
    U_\p \\ Y_\p \\ U_\f \\ Y_\f 
\end{bmatrix} g_1
= \begin{bmatrix}
    U_\p \\ Y_\p \\ U_\f \\ M + Y_\f(I-\Pi_1)
\end{bmatrix} g_1 
= \begin{bmatrix}
    U_\p \\ Y_\p \\ U_\f \\ M
\end{bmatrix} g_1
= \begin{bmatrix}
    u_\ini \\ y_\ini+\sigma_y \\ u \\ y
\end{bmatrix}, 
\]
where we have applied the fact that $(I-\Pi_1)g_1=0$ from \cref{eqn:equival-data-driven-SPC-2}. Thus, the set of variables $g_1, u_1, y_1$ and $\sigma_{y_1}$ is also a feasible solution for~$\cref{eqn:data-driven-spc}$.

We next assume $u_1, y_1, \sigma_{y_1}$ and $g_1$ is a feasible solution for \cref{eqn:data-driven-spc}. We define $\tilde{g}_1 = H_1^\dag 
\col(u_\ini, y_\ini + \sigma_{y_1}, u_1)$, which satisfies $y_1 = Y_\f \tilde{g}_1$ from \cref{eqn:relate-y-YF}. We first verify that $\tilde{g}_1$, together with $u_1, y_1, \sigma_{y_1}$, satisfies \cref{eqn:equival-data-driven-SPC-1}: 
\[
\begin{bmatrix}
    U_\p \\ Y_\p \\ U_\f \\ Y_\f
\end{bmatrix} \tilde{g}_1 = 
\begin{bmatrix}
    H_1   \\
    Y_\f
\end{bmatrix}
H_1^\dag \begin{bmatrix}
        u_\ini \\ y_\ini + \sigma_{y_1} \\ u_1
        \end{bmatrix} = \begin{bmatrix}
            u_\ini \\ y_\ini + \sigma_{y_1} \\ u_1 \\ y_1
        \end{bmatrix}.
\]
For the satisfaction of \cref{eqn:equival-data-driven-SPC-2}, since $\tilde{g}_1$ is in the range space of $H_1^\dag$ and $\Pi_1$ is the orthogonal projector onto the row space of $H_1$, we have $\Pi_1 \tilde{g}_1 = \tilde{g}_1$ (the range space of $H_1^\dag$ and row space of $H_1$ are equivalent; see \cref{eqn:p-range}), which implies $\|(I-\Pi_1)\tilde{g}_1\|_2 = \|\tilde{g}_1 -\tilde{g}_1\|_2 = 0$. Thus, $u_1, y_1, \sigma_{y_1}$ and $\tilde{g}_1$ is a feasible solution for~\cref{eqn:equival-data-driven-SPC}.

The uniqueness of the optimal solution $u^*, y^*$ and $\sigma_y^*$ for \cref{eqn:data-driven-spc} basically comes from strong convexity. For notational simplicity, we define $x= \col(u, y, \sigma_{y}, g)$ and $f(x) = \|u\|_R^2 + \|y\|_Q^2 + \lambda_y \|\sigma_y\|_2^2$ as the decision variable and objective function of \cref{eqn:data-driven-spc}, respectively. Suppose that $x_1$ and $x_2$ are two optimal solutions with different $u, y$ and $\sigma_y$. Let the optimal value be $f^*$. We then construct a convex combination $x_3 = \alpha x_1 + (1- \alpha) x_2$ where $0 < \alpha < 1$. This new point $x_3$ is also feasible since
\[
\begin{bmatrix}
    U_\p \\ Y_\p \\ U_\f \\ M
\end{bmatrix} g_3 =
\begin{bmatrix}
    U_\p \\ Y_\p \\ U_\f \\ M
\end{bmatrix} (\alpha g_1 + (1-\alpha)g_2) = 
\begin{bmatrix}
    u_\ini \\ y_\ini + \alpha \sigma_{y_1} + (1-\alpha) \sigma_{y_2} \\ \alpha u_1 + (1-\alpha)u_2 \\ \alpha y_1 + (1-\alpha)y_2 
\end{bmatrix} = 
\begin{bmatrix}
    u_\ini \\ y_\ini + \sigma_{y_3} \\ u_3 \\ y_3   
\end{bmatrix}. 
\] 
It is obvious that $f(x)$ is a strongly convex function with respect to $u, y$ and $\sigma_y$ and its value is not affected by $g$. Thus, 
we have
$
f(x_3) = f(\alpha x_1 + (1-\alpha) x_2) < \alpha f(x_1) + (1-\alpha) f(x_2) = f^*,
$ 
which contradicts our assumption. The optimal solution to \cref{eqn:data-driven-spc} is thus unique. This completes our proof.

\subsection{Proof of \Cref{prop:reg-equil-condtion}} \label{appendix:proposition}

The proof idea is to eliminate the variables $g$ and $\bar{g}$ in \cref{eqn:DeePC-Hybrid-2-norm} and \cref{eqn:DeePC-SVD-2-norm}, and to show the resulting problems are the same. Specifically, we will show that 
\begin{enumerate}
    \item Under \cref{eqn:prop-1}, upon fixing the values of $u,y,$ and $\sigma_y$, we can analytically find $g$ and $\bar{g}$ that minimize the objective functions in \cref{eqn:DeePC-Hybrid-2-norm} and \cref{eqn:DeePC-SVD-2-norm}.
    \item Under \Cref{eqn:prop-2}, the objective functions in \cref{eqn:DeePC-Hybrid-2-norm} and \cref{eqn:DeePC-SVD-2-norm} are the same after replacing $g$ and $\bar{g}$ with $u,y$ and $\sigma_y$. 
\end{enumerate}
With these two properties and the fact that \cref{eqn:DeePC-Hybrid-2-norm} and \cref{eqn:DeePC-SVD-2-norm} are strictly convex in $u,y$ and $\sigma_y$,  
we conclude that they have the same unique optimal solution $u^*, y^*$ and $\sigma_y^*$.

Let us prove the property 1.  Recall that we denote $v = \col(u_\ini,y_\ini+\sigma_y,u,y)$. Given the values of $u, y, \sigma_y$, without loss of generality, we can represent $g$ as 
\begin{equation}
\label{eqn:g-reform}
    g = g_\textrm{p} + \hat{g}, 
\end{equation}
where $g_\textrm{p} = H^\dag v$ and $\hat{g}$ is a vector in $\textrm{null}(H)$. We then substitute \cref{eqn:g-reform} into the regularizer $\|Gg\|_2^2$ in \cref{eqn:DeePC-Hybrid-2-norm}, leading to 
\[
\begin{aligned}
\|Gg\|_2^2 & = \|G(g_\textrm{p}+\hat{g})\|_2^2 \\
&= \|G g_\textrm{p}\|_2^2 + \|G \hat{g}\|_2^2 + 2(g_\textrm{p}^\tr G^\tr G\hat{g}) \\
&= \|G g_\textrm{p}\|_2^2 + \|G \hat{g}\|_2^2 + 2(v^\tr (H^\dag)^\tr G^\tr G\hat{g}).
\end{aligned}
\]
Note that we have $\textrm{rowsp}((H^\dag)^\tr G^\tr G) \subseteq \textrm{rowsp}(H)$, which is because $\textrm{rowsp}((H^\dag)^\tr) = \textrm{rowsp}(H)$ from \cref{eqn:p-range} and the property \cref{eqn:prop-1} that $\textrm{rowsp}(H G^\tr G) \subseteq \textrm{rowsp}(H)$. Then we get 
$$(H^\dag)^\tr G^\tr G\hat{g} = 0$$  
thanks to $\hat{g} \in \textrm{null}(H)$. Thus, we can further derive 
\begin{equation} \label{eqn:G-solution-for-g}
\|Gg\|_2^2 = \|G g_\textrm{p}\|_2^2 + \|G \hat{g}\|_2^2. 
\end{equation}

From \cref{eqn:G-solution-for-g}, we know that 
\begin{subequations}
\begin{equation} \label{eq:g-solution}
    g = g_\textrm{p}= H^\dag v
\end{equation}
is an optimal solution that minimizes the objective function in \cref{eqn:DeePC-Hybrid-2-norm}. 
On the other hand, since $\bar{H}$ has linearly independent columns in \cref{eqn:DeePC-SVD-2-norm}, $\bar{g}$ has a unique solution for fixed values of $u, y$ and $\sigma_y$. This unique solution is given by 
\begin{equation}\label{eq:gbar-solution}
\bar{g}_\textrm{p} = \bar{H}^\dag v.
\end{equation}
\end{subequations}

Next, we substitute \cref{eq:g-solution,eq:gbar-solution} into the objective functions in \cref{eqn:DeePC-Hybrid-2-norm} and \cref{eqn:DeePC-SVD-2-norm}, leading to  
\[
\begin{aligned}
f_1(u, y, \sigma_y)  
= & \ \|u\|_R^2 + \|y\|_Q^2 + \lambda_g \|G V \Sigma^{-1}W^\tr v\|_2^2 + \lambda_y \|\sigma_y\|_2^2, \\
f_2(u, y, \sigma_y) 
= & \ \|u\|_R^2 + \|y\|_Q^2 + \lambda_g \|\bar{G} \Sigma^{-1}W^\tr v\|_2^2 + \lambda_y \|\sigma_y\|_2^2.
\end{aligned}
\]
Under \cref{eqn:prop-2} that $V^\tr G^\tr G V  = \bar{G}^\tr \bar{G}$, we have
$
\|G V \Sigma^{-1}W^\tr v\|_2^2 = \|\bar{G} \Sigma^{-1}W^\tr v\|_2^2, \, \forall v.
$
This means that \cref{eqn:DeePC-Hybrid-2-norm} and \cref{eqn:DeePC-SVD-2-norm} have the same objective function in terms of $u, y, \sigma_y$. This completes our proof.

\section{Additional numerical results} \label{appendix:numerical-results}
In this section, we present numerical experiments with nonlinear systems that exhibit ``bias" errors. We conduct a numerical comparison between the indirect data-driven method and various variants of \DeePC{} across systems characterized by varying degrees of nonlinearity.

\noindent \textbf{Experiment setup.} We consider the nonlinear Lotka-Volterra dynamics used in \cite{dorfler2022bridging}
\begin{equation}
\label{eqn:nonlinear-sys}
\dot{x} = 
    \begin{bmatrix}
        \dot{x}_1 \\
        \dot{x}_2
    \end{bmatrix} = \begin{bmatrix}
        a x_1 - b x_1 x_2 \\
        d x_1 x_2 - c x_2 + u
    \end{bmatrix},
\end{equation}
where $x_1, x_2$ denote prey and predator populations and $u$ is the control input. We used $a =c=0.5, b = 0.025, d=0.005$ in our experiments. 

We first linearize the system \cref{eqn:nonlinear-sys} around the equilibrium $(\bar{u}, \bar{x}_1, \bar{x_2}) = (0, \frac{c}{d}, \frac{a}{b})$, which yields a linear system in the error state space. After discretization, we obtain a linear system as  
\[
\hat{x}(k+1) = f_{\text{linear}}(\hat{x}(k), \hat{u}(k)) = \begin{bmatrix}
    \hat{x}_1(k) + \Delta t (-b \bar{x}_1 \hat{x}_2(k)) \\
    \hat{x}_2(k) + \Delta t (d \bar{x}_2 \hat{x}_1(k) + \hat{u}(k))
\end{bmatrix}
\]
where $\Delta t = 0.1$ is the time step for discretization. We also discretize the nonlinear system in the error state space
\[
\begin{aligned}
\hat{x}(k+1) &= f_{\text{nonlinear}}(\hat{x}(k), \hat{u}(k)) \\
&= \begin{bmatrix}
    \hat{x}_1(k) + \Delta t (a(\hat{x}_1(k)+\bar{x_1})-b(\hat{x}_1(k)+\bar{x_1})(\hat{x}_2(k)+\bar{x_2})) \\
    \hat{x}_2(k) + \Delta t (d (\hat{x}_1(k)+\bar{x_1})(\hat{x}_2(k)+\bar{x_2})-c(\hat{x}_2(k)+\bar{x_2}) + \hat{u}(k))
\end{bmatrix}.
\end{aligned}
\]
We construct systems with various nonlinearity by interpolating between $f_{\textrm{linear}}$ and $f_{\textrm{nonlinear}}$ that is 
\[
\hat{x}(k+1) = \epsilon \cdot f_{\text{linear}}(\hat{x}(k), \hat{u}(k)) + (1-\epsilon) \cdot f_{\text{nonlinear}}(\hat{x}(k), \hat{u}(k)).
\]
The length of the pre-collected trajectory is $T=300$ and the prediction horizon, initial sequence are set to $N = 60$ and $T_\ini = 4$, respectively. We choose $Q=I, R=0.5I$ and $\mathcal{U} = [-20, 20]$. We further fix $\lambda_1 = 300, \lambda_2 = 100$ and $\lambda_y = 10000$ for all simulations. 

\vspace{2mm}

\begin{wrapfigure}[14]{r}{0.45\textwidth}
\vspace{-6mm}
    \includegraphics[width=0.45\textwidth]{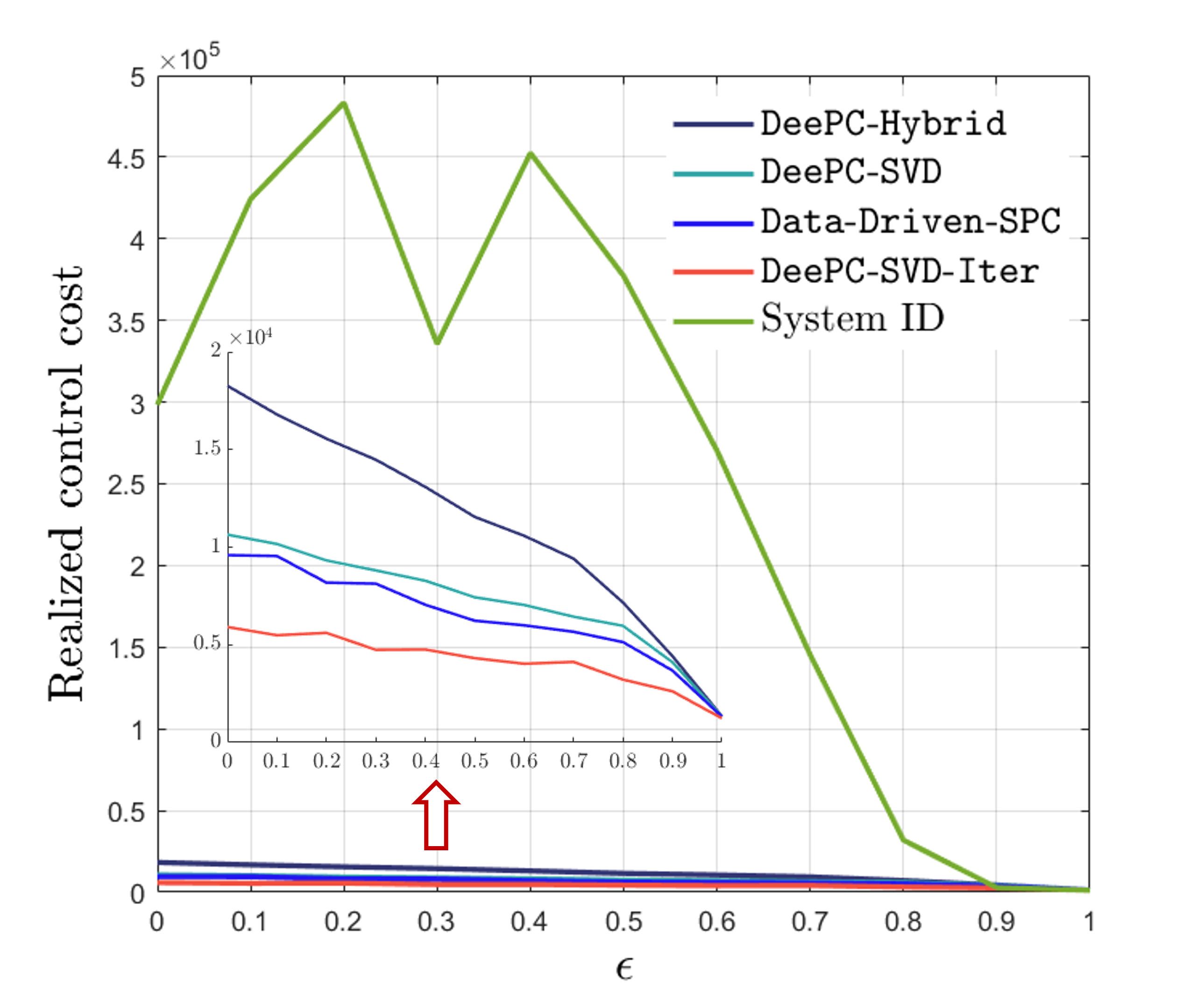}
    \vspace{-8mm}
    \caption{Comparison of realized cost for system ID and \DeePC{} based Methods}
    \label{fig:cost-nonlinear-All}
\end{wrapfigure}
\noindent \textbf{Comparison of direct/indirect methods.} We compare the realized control costs for the indirect system ID approach and different convex approximations on systems with varying degrees of nonlinearity. We use MATLAB function N4SID \citep{van1994n4sid} to identify the system and compute the optimal control input by solving \eqref{eqn:predictive-control}. Model orders are chosen to be 2 and 4 for \method{DeePC-SVD-Iter} in \Cref{alg:iter-SLRA} and N4SID, respectively, as they provide consistent good performance for all experiments.    

Similar to \Cref{sec:results}, we average the realized control costs over 100 pre-collected trajectories. We note that the identified model from N4SID is often ill-conditioned when the nonlinearity is high which caused numerical issues in solving \cref{eqn:predictive-control} in some of our experiments. We discard these infeasible solutions when computing the average performance for the indirect System ID approach. 

The results are shown in \Cref{fig:cost-nonlinear-All}. Both direct (convex approximations) and indirect (system ID) approaches perform well when the nonlinearity is low ($\epsilon \in [0.9, 1]$). However, the cost for the indirect method significantly increases with higher nonlinearity, while the performance of direct methods remains relatively consistent. Our numerical experiments suggest that the indirect system ID approach is inferior to direct data-driven methods in the case of ``bais'' error. The superior performance of direct data-driven methods is consistent with experimental observations from \cite{dorfler2022bridging}. The indirect system ID method projects the noisy data on a fixed linear model (which might induce ``bias'' error due to selecting a wrong model class); on the other hand, the complexity of the LTI system is regularized but not specified in direct methods, which provides more flexibility and leads to the superior performance for controlling nonlinear systems. Another intuition for this observation is that a nonlinear system in a finite-time horizon can in principle be approximated very well by an LTI system with sufficiently large dimension, which might be supported by lifting arguments from Koopman theory. 

\begin{figure}
    \centering
    \setlength{\abovedisplayskip}{0pt}
\includegraphics[width=0.78\textwidth]{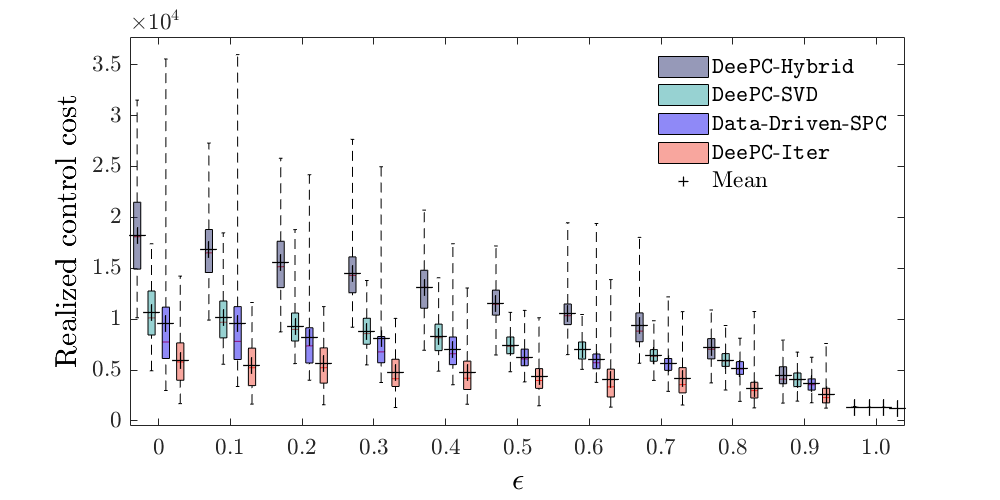}
    \caption{Comparison of realized control cost of different convex approximations for systems with varying nonlinearity. The navy block and cyan block represent \method{DeePC-Hybrid} and \DeePCsvd{}, respectively. The blue block and red block correspond to \method{Data-Driven-SPC} and \method{DeePC-SVD-Iter}.}
    \label{fig:cost-nonlinear}
\end{figure}

\vspace{2mm}
\noindent \textbf{Comparison of \DeePC{} variants.} We then compare the performance of different convex approximations as displayed in \Cref{fig:cost-nonlinear}. The \method{DeePC-SVD-Iter} and \method{Data-Driven-SPC} with structured data-driven representation outperform \method{DeePC-Hybrid} and \DeePCsvd{} which lift constraints as regularizers in objective functions. Furthermore, \method{DeePC-SVD-Iter} degrades much less than other methods. These numerical results suggest that we might obtain the additional benefits when employing appropriate techniques from system ID to pre-process the trajectory library of nonlinear systems. Our current results and analysis focus on open-loop control performance, and it will be very interesting to further investigate the performance of various convex approximations in the receding horizon closed-loop implementation. 

\end{document}